\RequirePackage{varioref}

\documentclass{siamart250211}

\usepackage{repdefs_siam}
\usepackage{subdepth}
\usepackage{amsmath}
\usepackage{amssymb}
\usepackage{mathtools}
\usepackage{longtable}
\usepackage{arydshln}
\usepackage{appendix}
\usepackage{color}
\usepackage[normalem]{ulem}

\newcommand{\xs}{x(\sigma)}
\newcommand{\xsi}{x(\sigma_i)}
\newcommand{\ys}{y(\sigma)}
\newcommand{\ysj}{y(\sigma_j)}
\newcommand{\sigmamin}{\sigma_{\min}}
\renewcommand{\sigmamax}{\sigma_{\max}}
\newcommand{\sigmas}{\sigma\sub{S}}

\newcommand{\rank}{\mbox{rank}}
\newcommand{\rankeps}{\rank_{\epsilon}}
\newcommand{\Sbar}{\bar{S}}
\newcommand{\As}{A_{\sigmas}}

\newcommand{\xt}{x_{\sigmas}(\theta)}

\newcommand{\PP}{\revised{R}}
\newcommand{\xp}{x\sub{\PP}}
\newcommand{\xps}{x\sub{\PP}(\sigma)}
\newcommand{\xpstar}{x\sub{\PP}{}_*}
\newcommand{\Ap}{A\sub{\PP}}
\newcommand{\bp}{b\sub{\PP}}
\newcommand{\sigmapmin}{\sigma\sub{\PP}{}_{\min}}
\newcommand{\kmax}{m}
\newcommand{\revised}[1]{{{#1}}}

\let\oldnl\nl% Store \nl in \oldnl
\newcommand{\nonl}{\renewcommand{\nl}{\let\nl\oldnl}}% Remove line # for 1 line

\def\arrvline{\hfil\kern\arraycolsep\vline\kern-\arraycolsep\hfilneg}

\usepackage[algo2e,linesnumbered,ruled,vlined]{algorithm2e}
\usepackage{setspace}
\SetKwInput{KwInput}{Input}
\SetKwInput{KwOutput}{Output}

\renewcommand{\thealgocf}{\arabic{section}.\arabic{algocf}}

\newcommand{\papertitle}{Extended-Krylov-subspace methods for
trust-region and norm-regularization subproblems}
\newcommand{\paperauthor}{Hussam Al Daas and Nicholas I. M. Gould\thanks{
   Computational Mathematics Theme,
   Scientific Computing Department,
   Rutherford Appleton Laboratory, 
   Chilton, OX11 0QX, England (UK).
   Email: hussam.al-daas@stfc.ac.uk, nick.gould@stfc.ac.uk}}

\title{\papertitle}
\author{\paperauthor}

\begin{document}

\maketitle
%\begin{center} Preprint STFC-P-2025-002 \end{center}
\begin{abstract}
\noindent
We consider an effective new method for solving trust-region and 
norm-regularization problems that arise as subproblems in many optimization
applications. We show that the solutions to such subproblems
\revised{effectively lie in a very-low-dimensional subspace}
%\sout{on a manifold of approximately very low rank}
as a function of their
controlling parameters (trust-region radius or regularization weight).
Based on this, we build a basis \revised{spanning these solutions}
%\sout{for this manifold}
using an efficient
extended-Krylov-subspace iteration that involves a single matrix 
factorization. The problems within the subspace using such a basis
may be solved at very low cost using effective high-order root-finding methods.
This then provides an alternative to common methods using
multiple factorizations or standard Krylov subspaces. We provide numerical
results to illustrate the effectiveness of our {\tt TREK}/{\tt NREK} 
approach.
\end{abstract}

%\vspace*{-15mm}
\numsection{Introduction}

Given an $n$ by $n$  symmetric (Hessian) matrix $A$, an $n$-vector $b$, a 
scalar ``radius'' $\Delta > 0$, and the Euclidean norm $\|\cdot\|$, 
the trust-region subproblem is to find
\eqn{tr}{x_* = 
 \argmininx{} \half x^T A x - b^T x \tim{subject to} \|x\| \leq \Delta.}
Such problems arise in trust-region methods for (non-convex) nonlinear 
optimization \cite{ConnGoulToin00,NoceWrig99}. 
A variety of methods have been proposed to solve them, some based on
matrix factorization \cite{GoulRobiThor10:mpc,MoreSore83:sisc}, others
on approximations from Krylov subspaces 
\cite{CartGoulToin22,GoulLuciRomaToin99:siopt}
and some on reformulation as eigenproblems 
\cite{AdacIwatNakaTake17:siopt,RendWolk97:mp,RojaSantSore01:siopt,SantSore95}.
It is well known \cite{Gay81:sisc,MoreSore83:sisc} that the solution $\xs$ 
satisfies the first-order optimality condition
\eqn{fon}{ (A + \sigma I) \xs = b,}
where $\|\xs\| \leq \Delta$ and the scalar shift
\eqn{sigmamin}{\sigma \geq \sigmamin \eqdef \max(0,- \lambda_1(A))}
satisfies the complementarity condition $\sigma (\|\xs\| - \Delta) = 0$; 
here $\{\lambda_1(A), \ldots, \lambda_n(A)\}$ are the eigenvalues of $A$, 
in increasing order, i.e., $\lambda_i(A) \leq \lambda_{i+1}(A)$, for
$i = 1,\ldots,n-1$. That is, either, for positive semi-definite $A$, 
\eqn{trivial}{ A x(0) = b \tim{with} \|x(0)\| \leq \Delta}
or, with unrestricted $A$,
\eqn{ls}{ ( A + \sigma I ) \xs = b \tim{with} \|\xs\| = \Delta
\tim{for some} \sigma \geq \sigmamin.}
We refer to the first possibility as the trivial---or interior---case, while
the second is the general one; there is an additional probability-zero
possibility, commonly called the hard case---this occurs if $b$ is orthogonal
to all the eigenvectors for the leftmost eigenvalue $\lambda_1(A)$ 
and the radius is too large---that we shall ignore for the time 
being\footnote{A few methods are specifically designed to cope with this
possibility \cite{BeckVasi18:siopt,CarmDuch20:sirev}.}.
If we can rule out the trivial and hard cases, \req{ls}
then suggests an iteration in which, for a sequence of positive $\sigma$,
$\xs$ is found by solving the linear system \req{fon}, and its suitability
assessed by examining the size of the residual $ \|\xs\| - \Delta$.

An interesting question is, therefore, what we can say about the manifold 
$\calM(A,b) = \{ x \st x = \xs \tim{for some} 
 \sigma \in [\sigmamin,\infty)\}$
or more specifically the \revised{subspace} 
%\sout{manifold}
$\calS(A,b) = \mbox{span}\{\xs\}$, that contains $\calM(A,b)$, 
as $\sigma$ varies between $\sigmamin$ and $+\infty$.
There is no loss of generality, from a theoretical perspective, 
in considering a problem with a diagonal Hessian, since if $A$ has 
the spectral factorization 
\eqn{spectral}{A = U \Lambda U^T,} 
involving an orthogonal matrix $U$ whose 
columns are its eigenvectors, and a diagonal matrix $\Lambda$ whose entries 
are its eigenvalues, \req{tr} is equivalent to
\eqn{trt}{\revised{\argminin{y\in \tinyRe^n}}
 \half y^T \Lambda y - d^T y \tim{subject to}  \|y\| \leq \Delta}
under the transformation 
\disp{y = U^T x, \tim{where} d = U^Tb.}
%\vspace*{-8mm}

Consider, then, problem \req{trt} with $n = 10^4$, $d$ a vector of ones,
eigenvalues $\lambda_i$, $i = 1,\ldots,n$,  in $(0,1]$ that are either
(i) equi-spaced at $i / n$, 
(ii) clustered at the lower end at $i^2/n^2$, 
(iii) clustered at the higher end at $1 - (i-1)^2/n^2$, or
(iv) logarithmically-distributed in $(10^{-15},1]$ as
$\exp(\log(10^{-15}) + (i-1) \beta)$, where 
$\beta =  (\log(1) - \log(10^{-15})/(n - 1)$, and $p = 10n$ 
evaluation points $\sigma = \sigma_j = j/n$ for $j = 1,\ldots,p$ in $(0,10]$.
In this case, the components of $\ys$ from the first-order optimality system
\disp{ (\Lambda + \sigma I) \ys = d}
for \req{trt} are $\revised{[y(\sigma)]_i} = 1 / (\lambda_i  + \sigma)$. 
If we approximate the \revised{basis of the subspace}
%\sout{solution manifold}
$\calS(\Lambda,d)$ by the discrete set $\{ \ysj, j = 1,\ldots,p\}$
at the evaluation points, and compute the singular value decomposition of
the resulting $n$ by $p$ matrix $S$ with entries
\eqn{cauchy-matrix}{\revised{[S]_{i,j} = [d]_i} / (\lambda_i  + \sigma_j) 
 \tim{for} i = 1, \ldots, n \tim{and} j = 1, \ldots, p, }
we find 37 significant singular values (those that are
larger than $10^{-15}$ of the largest) in case (i), and for cases (ii)--(iv)
there are 38, 36 and 35 respectively. The precise singular values are given
in Table~\ref{sv-table}.

\begin{table}[htbp]
\begin{center}
{\footnotesize
\begin{tabular}{|cc|cc|cc|cc|}
\hline
\multicolumn{2}{|c}{(i): equi-spaced} & 
\multicolumn{2}{|c}{(ii): clustered lower} & 
\multicolumn{2}{|c|}{(iii): clustered higher} &
\multicolumn{2}{|c|}{(iv): logarithmically} \\
%\mbox{(i)} & \mbox{(i)} & \mbox{(i)} \\
\hline
  2.53e+04 & 1.32e-03  & 1.34e+05 & 6.33e-03 & 1.86e+04 & 4.84e-04 & 
  1.09e+06 & 4.64e-03 \\ 
  1.44e+04 & 4.99e-04  & 4.38e+04 & 2.48e-03 & 1.09e+04 & 1.74e-04 & 
  7.32e+04 & 1.81e-03 \\
  6.89e+03 & 1.87e-04  & 1.92e+04 & 9.70e-04 & 5.12e+03 & 6.23e-05 & 
  2.09e+04 & 7.04e-04 \\
  3.07e+03 & 7.00e-05  & 8.81e+03 & 3.78e-04 & 2.21e+03 & 2.22e-05 & 
  7.75e+03 & 2.72e-04 \\
  1.33e+03 & 2.60e-05  & 3.95e+03 & 1.46e-04 & 9.19e+02 & 7.85e-06 & 
  3.23e+03 & 1.05e-04 \\
  5.64e+02 & 9.64e-06  & 1.73e+03 & 5.65e-05 & 3.76e+02 & 2.77e-06 & 
  1.39e+03 & 4.04e-05 \\
  2.36e+02 & 3.55e-06  & 7.42e+02 & 2.17e-05 & 1.51e+02 & 9.70e-07 & 
  5.96e+02 & 1.55e-05 \\
  9.78e+01 & 1.30e-06  & 3.15e+02 & 8.33e-06 & 6.03e+01 & 3.39e-07 & 
  2.52e+02 & 5.90e-06 \\
  4.01e+01 & 4.76e-07  & 1.33e+02 & 3.18e-06 & 2.38e+01 & 1.18e-07 & 
  1.05e+02 & 2.24e-06 \\
  1.63e+01 & 1.73e-07  & 5.54e+01 & 1.21e-06 & 9.27e+00 & 4.07e-08 & 
  4.35e+01 & 8.51e-07 \\
  6.55e+00 & 6.29e-08  & 2.29e+01 & 4.59e-07 & 3.58e+00 & 1.40e-08 & 
  1.79e+01 & 3.21e-07 \\
  2.61e+00 & 2.27e-08  & 9.44e+00 & 1.74e-07 & 1.37e+00 & 4.80e-09 & 
  7.30e+00 & 1.21e-07 \\
  1.04e+00 & 8.18e-09  & 3.86e+00 & 6.54e-08 & 5.20e-01 & 1.64e-09 & 
  2.96e+00 & 4.54e-08 \\
  4.07e-01 & 2.93e-09  & 1.57e+00 & 2.46e-08 & 1.96e-01 & 5.57e-10 & 
  1.19e+00 & 1.70e-08 \\
  1.59e-01 & 1.05e-09  & 6.34e-01 & 9.19e-09 & 7.33e-02 & 1.89e-10 & 
  4.79e-01 & 6.35e-09 \\
  6.18e-02 & 3.73e-10  & 2.55e-01 & 3.43e-09 & 2.72e-02 & 6.37e-11 & 
  1.91e-01 & 2.36e-09 \\
  2.38e-02 & 1.32e-10  & 1.02e-01 & 1.28e-09 & 1.00e-02 & 2.14e-11 & 
  7.60e-02 &         \\
  9.14e-03 & 4.69e-11  & 4.05e-02 & 4.74e-10 & 3.67e-03 &          & 
  3.01e-02 &          \\
  3.48e-03 &           & 1.61e-02 & 1.76e-10 & 1.34e-03 &          & 
  1.18e-02 &          \\
\hline
\end{tabular}
}
\caption{\revised{Dominant} singular values of 
\revised{$S$ in~\req{cauchy-matrix}}
%\sout{dominant manifold subspace} 
for four eigenvalue distributions}
\label{sv-table}
\end{center}
\end{table}
It is then remarkable that \revised{$\calS(\Lambda,d)$} 
%\sout{the manifold} 
appears to be spanned by so few
significant vectors, that is, it is approximately of very low rank.
This turns out not to be a coincidence, and we shall exploit this 
throughout.

This paper is organised as follows. 
In Section~\ref{low-rank-section}, we explain why the 
\revised{subspace spanned by the solutions} 
%\sout{manifold} 
of interest is (approximately) low \revised{dimensional.} 
%\sout{rank}
We follow this in Section~\ref{manifoild-approximation-section}
by considering methods that build approximations to this \revised{subspace}.
%\sout{manifold}
A particularly appealing approach may be constructed from a so-called 
extended Krylov subspace, and we show how in Section~\ref{EKS-section}.
This then provides the basis for methods aimed at the trust-region subproblem
\req{tr}, and we provide details for positive-definite $A$ in 
Section~\ref{definite-subproblem-solve-section}. We follow this with
a description of how we may extend this for general $A$ in
Section~\ref{general-subproblem-solve-section}. 
Numerical experiments follow in Section~\ref{experiments-section}.
We show that the same ideas trivially extend to the regularized-norm
problem in Section~\ref{regularization-section},
and finally we conclude in
Section~\ref{conclusions-section}.

\numsection{Cauchy matrices and accurate low-rank approximation}
\label{low-rank-section}

The matrix \revised{$S$ in~}\req{cauchy-matrix} is an example of what is commonly
known as a generalised Cauchy---or Cauchy-like---matrix \cite{Cauc24:wiki}.
To be more specific, suppose that we know the ordered, increasing eigenvalues
$\lambda_i \in [\lambda_1, \lambda_n]$ of $A$ (and thus $\Lambda$), 
and have picked a set of distinct shifts
$\sigma_i \in [\sigmamin^+, \sigmamax]$, $i = 1, \ldots, p$, where
\disp{\sigmamin^+ \eqdef \max(0,- \lambda_1+ \delta) \tim{and}
\sigmamax > \sigmamin^+}
for some $\delta > 0$. Then
$\sigmamin = \max(0,- \lambda_1) \leq \sigmamin^+$, and moreover
$\sigmamin^+ \geq 0 > - \delta \geq -\lambda_1$ whenever 
$\lambda_1 \geq \delta$ and 
$\sigmamin^+ \geq - \lambda_1+ \delta > - \lambda_1$ otherwise.
Thus $-\sigmamin^+ < \lambda_1$, indeed $\sigmamin^+ + \lambda_1 \geq \delta$, 
and the intervals 
$[-\sigmamax, -\sigmamin^+]$ and $[\lambda_1, \lambda_n]$ are non-overlapping.

%\revised{NICK: MAYBE MY ENGLISH IS NOT HELPING MUCH HERE, 
%BUT I DON'T FEEL THE FOLLWOING SENTENCE IS RIGHT} 
\revised{
Clearly rows $i$ and \revised{$\ell$} of $S$ are linearly dependent if 
$\lambda_i = \revised{\lambda_{\ell}}$, and thus  $\rank(S) = \rank(\Sbar)$,
where $\Sbar$ comprises the rows of $S$ in which each set of such dependencies
is replaced by a single representative.
To investigate the rank of $S$, we thus restrict our attention to $\Sbar$.}
For any $m$ by $n$ matrix $X$ with decreasing singular values 
$\sigma_1(\revised{X}) \geq \sigma_2(X) \geq \ldots \geq \sigma_k(X)$ 
for increasing $k$,
the $\epsilon$-rank, $\rankeps(X)$, of $X$ is the smallest integer $k$ for which
$\sigma_{k+1}(X) \leq \epsilon \sigma_1(X)$.
Then Beckermann and Townsend 
\cite[\S4.2]{BeckTown17:simaa} show that\footnote{Stronger bounds
are possible \cite{BeckTown17:simaa}
in terms of the so-called Gr\"{o}tzsch ring function, and other bounds
are also known, e.g. \cite{Gras04:nlaa}.}
\disp{\rankeps(\Sbar) \leq \left\lceil
\log(16\gamma) \log(4/\epsilon) / \pi^2
\right\rceil,}
involving the cross-ratio
\disp{\gamma = \left| \frac{(\sigmamax+\lambda_1)(\sigmamin^++\lambda_n)}
{(\sigmamax+\lambda_n)(\sigmamin^++\lambda_1)}\right|,}
because we have removed repeated eigenvalues, the shifts are distinct and 
the two sets are non-overlapping. 
%{\em Should we state this as a theorem, so that others can quote it?}
Notice that the cross ratio is bounded
from above as $\sigmamin^+ + \lambda_1 \geq \delta > 0$, and approaches 
$(\sigmamin^++\lambda_n)/(\sigmamin^++\lambda_1)$ as 
$\sigmamax \rightarrow \infty$; for positive definite $A$, the latter is
simply the spectral condition number. 
In particular, the log dependence on both the cross ratio and the precision
$\epsilon$ are significant in keeping the bound small, and
its dimension independence is also noteworthy.
For the examples illustrated in Table~\ref{sv-table}, this bound is 44, 
a very reasonable (upper-bound) approximation of the computed values.
%log(160000) * log(4/1.0e-15)/pi^2$

\numsection{\revised{Subspace} % \sout{Manifold}
approximation from Krylov subspaces}
\label{manifoild-approximation-section}

Now that we know that the solution we seek to \req{tr} {\em for any} 
$\Delta$ lies in an approximately low-dimensional subspace, our task is to
identify a suitable (for convenience, orthogonal) $n$ by $\ell$ basis 
matrix $V_{\ell}$ for this subspace. Once we do so, we may recover an
approximate solution $x_{\ell} = V_{\ell} y_{\ell}$, where
\eqn{trs}{y_{\ell} =  \revised{\argminin{y \in \tinyRe^{\ell}}}
 \half y^T P_{\ell} y - b_{\ell}^T y  
 \tim{subject to} \|y\| \leq \Delta\revised{,}}
with $P_{\ell} = V_{\ell}^T A V_{\ell}$ and $b_{\ell} = V_{\ell}^T b$. 
The fact that we expect $\ell$ to be small provides an opportunity to 
use methods that would be inappropriate has it been larger.

An obvious approach might be to simply select a set of 
$\sigma_i \geq \sigmamin^+$ for $i=1,\ldots,\ell$, 
compute the resulting solutions $\xsi$
and build an orthogonal basis of the set $\{\xsi\}$ by, for example,
the modified Gram-Schmidt process.
This has two fundamental issues. Firstly, it requires that we 
solve $\ell$ shifted linear systems, and this is not without cost.
Secondly, it is not known {\em a priori} how big 
$\ell$ should be. The latter may be addressed as follows.
Suppose that $y_{\ell}$ is the solution to \req{trs} with
$\|y_{\ell}\| = \Delta$. If so, there is a shift $\sigma_{\ell} \geq 0$, 
just as in \req{fon}, for which
\eqn{errory}{0 = (P_{\ell} + \sigma_{\ell}I) y_{\ell} - b_{\ell}
= V_{\ell}^T \left[ ( A + \sigma_{\ell}I)x_{\ell} - b\right]}
where $x_{\ell} = V_{\ell} y_{\ell}$ and $\|x_{\ell}\| = \Delta$,
since $V_{\ell}$ has orthogonal columns (i.e., $V_{\ell}^T V_{\ell} = I$).
As the solution we seek satisfies $( A + \sigma I)x = b$,
and as $V_{\ell}$ will be an increasingly better approximation to the 
low-dimensional subspace that contains the required solution, it suffices to
monitor the residual $r_{\ell} = ( A + \sigma_{\ell}I)x_{\ell} - b$ and to stop
when $\|r_{\ell}\|$ is sufficiently small.
This, however, does not avoid the cost of multiple factorizations; rival
methods, e.g., \cite{GoulRobiThor10:mpc,MoreSore83:sisc}, also use 
factorizations, but experience suggests that relatively few often suffice.
Another approach for constructing a good approximation of 
\revised{$\calS(A,b)$} 
%\sout{the manifold}} 
is to use a method based on a rational Krylov subspace
\cite{Ruhe84:LAA,AldaPali25:ARXIV}. However, this will still require
solving $\ell$ shifted linear systems.

The alternative is to try to build a good approximation to 
\revised{$\calS(A,b)$} 
%\sout{the manifold}
by other means. We shall resort to methods based on the extended Krylov subspace
(see Section~\ref{EKS-section}).
%Why not simply use the inverse subspace $\{ b, A^{-1}b, A^{-2}b, \ldots \}$?

To justify this, first note that given a shift $\sigma$, 
\req{fon} and \req{spectral} give that
\disp{x(\sigma) = (A+\sigma I)^{-1}b = U y(\sigma),}
where
\disp{ \revised{[y(\sigma)]_i} 
 = \frac{\revised{[d]_i}}{\lambda_i + \sigma},\;\; i=1,\ldots,n,  
\tim{and} d=U^Tb.}
Consider two cases and first suppose that 
\revised{$|\sigma| > |\lambda_i|$ for all $i=1,\ldots,n$.} 
Then using a Laurent series \cite{Laur24:wiki}, we have 
\revised{
\disp{[y(\sigma)]_i = [d]_i \sum_{j=1}^\infty
 \left(\frac{-\lambda_i}{\sigma}\right)^j.}
}
The smaller the ratio \revised{$|\lambda_i/\sigma|$} the faster the sequence
\revised{$(\left(-\lambda_i / \sigma \right)^j)_{j\in \mathbb{N}}$}
reaches machine precision and few
terms would be sufficient to approximate $y_i(\sigma)$. In other
words, $y(\sigma)$ can be very well approximated using a polynomial
with a moderate degree in \revised{$\Lambda$}, and hence $x(\sigma)$ can be well
approximated using a polynomial with a moderate degree in \revised{$A$}
applied to $b$.  
\revised{At the other extreme}, when \revised{$|\sigma|<|\lambda_i|$ for all} 
$i=1,\ldots,n$ \revised{and $A$ is non-singular}, a similar argument using \revised{the Taylor series 
\disp{[y(\sigma)]_i = [d]_i \sum_{j=0}^\infty
 \left(\frac{-\sigma}{\lambda_i}\right)^j.}
}
shows that $x(\sigma)$ can be well approximated
using a polynomial with moderate degree in \revised{$A^{-1}$} applied to $b$.

An extended Krylov subspace is by definition the space that constructs
the sum space of both Krylov subspaces 
$\calK(A,b) = \text{span}\{A^{i}b\}_{i \geq 0}$ and 
$\calK(A^{-1},b) = \text{span}\{A^{-i}b\}_{i \geq 0}$,
and constructs the aforementioned polynomials required to handle the two
ends of the behaviour spectrum.
Hence, these observations provide some intuition as to why we choose 
an extended Krylov method (see Section~\ref{EKS-section}) 
to solve the trust-region subproblem. 
Similar arguments have been made by Druskin and Knizhnerman 
\cite{DrusKniz98:simaa} for different applications.

For a better understanding of the convergence of the extended-Krylov-subspace 
method when solving shifted linear systems with a symmetric
positive definite matrix whose spectrum is contained in
$[\lambda_1,\lambda_n]$ and shifts lying outside
$(-\lambda_n-\delta_n,-\lambda_1+\delta_1)$ for $\delta_1,\delta_n > 0$,
we provide a convergence analysis in 
Appendix~\ref{sec:convergence analysis} that is heavily
based on that by Knizhnerman and Simoncini \cite{KnizSimo11:NUMMAT}.
This shows that there will
be an element $x_i(\sigma)$ of the extended Krylov subspace $\calK_{i}(A,b)$, 
defined in \req{eks}, for which 
\[\|x_i(\sigma) - \xs\| \leq c \mu^{i}\revised{(\sigma)}\]
for some constants $c$ and $\mu\revised{(\sigma)} < 1$ and all $i \geq 0$.

%\noindent
%{\em 
%Hussam: Is there anything that can be told about $\sigma+\lambda_1$ from the
%optimisation point of view? I believe when $\Delta$ get larger and
%larger, $\sigma_*$ approaches $\sigmamin$, that is, $\sigma_*+\lambda_1$
%approaches $\sigmamin$.
%Nick: That is correct.
%}

%In~\cite{KnizSimo11:NUMMAT}, the authors show that $B$ attains its maximum at $t = \left(\frac{\kappa^{1/4}-1}{\kappa^{1/4}+1}\right)$ with a value $\left(\frac{\kappa^{1/4}-1}{\kappa^{1/4}+1}\right)^2$, where $\kappa = \frac{\lambda_n}{\lambda_1}$.
%Therefore, 
%$$ \mu \leq \left(\frac{\kappa^{1/4}-1}{\kappa^{1/4}+1}\right)$$
%and
%$$\|x-x_i\| \leq C \left(\frac{\kappa^{1/4}-1}{\kappa^{1/4}+1}\right)^{2i}.$$
\numsection{The EKS algorithm}
\label{EKS-section}

Our aim is thus to construct efficiently an orthogonal basis  $V_{2m+1}$ 
of the evolving extended Krylov subspaces
\eqn{eks}{\calK_{2m}(A,b) = 
 \text{span}\{b,A^{-1}b,A b,A^{-2}b,A^2b,\ldots,A^{-m}b\}}
and
\disp{\calK_{2m+1}(A,b) = \revised{\text{span}\{\calK_{2m}(A,b),A^mb\}}.}
Fortunately, much of the groundwork is already in place. We suppose,
for this section and the next, 
that $A$ is positive definite\footnote{\revised{Formally,
Algorithm~\ref{algorithm-eks} merely requires that $A$ be non-singular, but
the matrices $P_{\ell}$ we subsequently develop are not always well defined
for indefinite $A$.
Fortunately, restricting attention to the definite case turns out to be 
insignificant,  since we show in Section~\ref{general-subproblem-solve-section}
how an inexpensive diagonal perturbation to $A$ allows us to solve 
both convex and nonconvex problems.}}. 
Then Jagels and Reichel\footnote{They also propose a second, more expensive 
algorithm that copes with indefinite $A$, \revised{and overcomes the defect 
mentioned in the previous footnote. Another variant that has a higher complexity and avoids the mentioned defect is presented in~\cite[Section 2 and 3]{Simo07:SISC}.}}
\cite{JageReic09:laa} have shown how to build such a basis 
if the order\footnote{We shall say why  our order is more convenient later.} 
in which the components $A^{-m}b$ and $A^mb$ arise is 
interchanged, that is if instead
\disp{\calK^{\prime}_{2m}(A,b) = 
 \text{span}\{b,A b,A^{-1} b,A^2b,A^{-2}b,\ldots,A^mb\}} 
and
\disp{\calK_{2m+1}(A,b) = \revised{\text{span}\{\calK^{\prime}_{2m}(A,b),A^{-m}b\}}.}
Using the same logic as in \cite[Alg.2.1]{JageReic09:laa},
that applies to $\calK^{\prime}_{2m}(A,b)$, it is
straightforward to derive the following algorithm, 
Algorithm~\ref{algorithm-eks} \vpageref{algorithm-eks},
that instead uses the order $\calK_{2m}(A,b)$; for completeness, 
we provide details of an alternative justification, that avoids Laurent 
polynomials, in Appendix~\ref{sec:orthogonal_basis}.
%\newpage

\begin{algorithm2e}[!ht]
\caption{The {\tt EKS} algorithm to find an orthogonal basis of 
$\calK_{2m+1}(A,b)$
\label{algorithm-eks}}
\KwInput{symmetric $A \in \Re^{n\times n}$, $b \in \Re^n$ and subspace 
dimension bound $m\geq 1$.}
\KwOutput{orthonormal $V_{2m+1} =$ \\ \hspace{11mm} 
   $( v_0 : v_{-1} : v_1 : \cdots : v_{k-1} : v_{-k} : v_{k} : v_{-k-1} : v_{k+1} : 
    \cdots : v_{-m} : v_m ) \in \Re^{n\times (2m+1)}.$ \\ {\bf Algorithm:}}
  \setstretch{0.8}
  \vspace{2mm}
  $\delta_0 \define \|b\|$\; 
  {\bf if} $\delta_0=0$ {\bf stop}, subspace completed\;
  $v_0 \define b/\delta_0$,
  $u \define A^{-1} v_0$\;
  $\beta_0 \define u^T v_0$, 
  $u \redef u - \beta_{0} v_{0}$\;
  $\delta_{-1} \define \|u\|$\;
  \For{$k = 1, 2, \ldots , m$}{
    {\bf if} $\delta_{-k}=0$ {\bf stop}, subspace completed\;
    $v_{-k} \define u/\delta_{-k}$,
    $u \define Av_{-k}$\;
    $\alpha_{k-1} \define u^T v_{k-1}$,
    $u \redef u - \alpha_{k-1} v_{k-1}$\;
    $\alpha_{-k} \define u^T v_{-k}$,
    $u \redef u - \alpha_{-k} v_{-k}$\;
    $\delta_{k} \define \|u\|$\;
    {\bf if} $\delta_{k}=0$ {\bf stop}, subspace completed\;
    $v_{k} \define u/\delta_{k}$,
    $u \define A^{-1} v_{k}$\;
    $\beta_{-k} \define u^T v_{-k}$,
    $u \redef u - \beta_{-k} v_{-k}$\;
    $\beta_{k} \define u^T v_{k}$,
    $u \redef u - \beta_{k} v_{k}$\;
    \If{$k<m$}{$\delta_{-k-1} \define \|u\|$\;
}}
\end{algorithm2e}
%\vspace*{-2mm}

\LinesNotNumbered

\noindent
Thus
\begin{align}
\delta_0 v_{0} & = b \nonumber\\
\delta_{-1} v_{-1} & = A^{-1} v_{0} - \beta_{0} v_{0} \label{it-zero} \\
\delta_{k} v_{k} & = Av_{-k} - \alpha_{-k} v_{-k} - \alpha_{k-1} v_{k-1} 
 \tim{for $k = 1, \ldots, \kmax$} \label{it-pos} \\
\delta_{-k-1} v_{-k-1} & = A^{-1} v_{k} - \beta_{-k} v_{-k} - \beta_{k} v_{k}
 \tim{for $k = 1, \ldots, \kmax-1$} \label{it-neg}
\end{align}
Each iteration of the main loop requires a solve and a product 
involving $A$, the initialization requires a single solve. The
remaining operations simply involve scalar products and vector sums.
The decision whether and when to stop in the main loop is the focus of
Sections~\ref{definite-subproblem-solve-section} and 
\ref{general-subproblem-solve-section}, but observe that the extended
Krylov subspace has definitely been completed if ever $\delta_i = 0$.

We consider the symmetric $\Re^{(2m+1)\times(2m+1)}$ ``projected'' matrix 
\[ P_{2m+1} = V_{2m+1}^T A V_{2m+1},\]
whose $i,j$-entry we denote as $p_{i,j} \revised{\equiv [P_{2m+1}]_{i,j}}$, 
as well as its leading  $2k$ and $2k+1$ sub-blocks
\[ P_{2k} = V_{2k}^T A V_{2k} \tim{and} P_{2k+1} = V_{2k+1}^T A V_{2k+1},\]
where
\[V_{2k} = ( v_0 : v_{-1} : v_1 : \cdots : v_{k-1} : v_{-k} )
 \tim{and} V_{2k+1} = ( V_{2k} : v_k),\]
that play the roles of \revised{$V_{\ell}$} in \req{trs}. 
It follows from \eqref{it-pos} that
\begin{equation}
Av_{-k} = \alpha_{k-1} v_{k-1} + \alpha_{-k} v_{-k} + \delta_{k} v_{k}
 \tim{for $k = 1, \ldots, \kmax$} \label{Amv}
\end{equation}
and the orthogonality of the $v_i$ then gives
\begin{equation}
p_{2k-1,2k} = v_{k-1}^T Av_{-k} = \alpha_{k-1}, \;\;
p_{2k,2k} = v_{-k}^T Av_{-k} = \alpha_{-k}, \;\;
p_{2k+1,2k} = v_{k}^T Av_{-k} = \delta_{k} \label{odd}
\end{equation}
and
\[
p_{i,2k} = v_i^T Av_{-k} = 0 \tim{for $i \neq 2k-1, \revised{2k}$ and $2k+1$.}
\]
By symmetry of $P_{2m+1}$, we also have
\[
p_{2k,2k-1} = \alpha_{k-1}, \;\;
p_{2k,2k} = \alpha_{-k}, \;\;
p_{2k,2k+1} = \delta_{k} \tim{and}
p_{2k,i} = 0 \tim{for the remaining $i$.}
\]
Thus, from \req{Amv} and \req{odd}, we have
\eqn{Amvk}{Av_{-k} = p_{2k,2k-1} v_{k-1} + p_{2k,2k} v_{-k} + p_{2k,2k+1} v_{k}
 \tim{for $k = 1, \ldots, \kmax$,}}
where formally any $p_{2k,i} = 0 $ if $i \leq 0$.
Similarly, \eqref{it-neg} gives
\[
A^{-1} v_{k} = \beta_{-k} v_{-k} + \beta_{k} v_{k} + \delta_{-k-1} v_{-k-1}
\tim{for $k = 1, \ldots, \kmax-1$} 
\]
and thus that
\eqn{bAvk}{
\arr{rl}{
\beta_{k} A v_{k} \ngap & = v_{k} - \beta_{-k} A v_{-k} - \delta_{-k-1} A v_{-k-1} \\
 & = 
- \beta_{-k} \alpha_{k-1} v_{k-1} - \beta_{-k} \alpha_{-k} v_{-k} 
   + (1  - \beta_{-k} \delta_{k} - \delta_{-k-1} \alpha_{k}) v_{k} \\
 & \gap\bgap 
- \delta_{-k-1} \alpha_{-k-1} v_{-k-1} - \delta_{-k-1} \delta_{k+1} v_{k+1}
}}
using \eqref{Amv}. Again, the orthogonality of the $v_i$ and symmetry 
implies that
%\footnote{When $\delta_{-k-1} = 0$, the last three terms in 
%\eqref{even} are $p_{2k+1,2k+1} = (1 - \beta_{-k} \delta_{k})/\beta_{k}$ and 
%$p_{2k+2,2k+1} = p_{2k+3,2k+1} = 0$, and thus do not require that we calculate
%$\alpha_{k}$, $\alpha_{-k-1}$ and $\delta_{k+1}$.}
\eqn{even}{\arr{rl}{
p_{2k-1,2k+1} \ngap & = p_{2k+1,2k-1} = v_{k-1}^T A v_{k}
 = - \bigfrac{\beta_{-k} \alpha_{k-1}}{\beta_{k}}, \\
p_{2k,2k+1} \ngap & = p_{2k+1,2k} 
 = v_{-k}^T A v_{k} = - \bigfrac{\beta_{-k} \alpha_{-k}}{\beta_{k}}, \\
p_{2k+1,2k+1} \ngap & = v_{k}^T A v_{k} 
 = \bigfrac{1 - \beta_{-k} \delta_{k} - \delta_{-k-1} \alpha_{k}}{\beta_{k}}, \\
p_{2k+2,2k+1} \ngap & = p_{2k+1,2k+2} = v_{-k-1}^T A v_{k} 
 = - \bigfrac{\delta_{-k-1} \alpha_{-k-1}}{\beta_{k}}, \\
p_{2k+3,2k+1} \ngap & = p_{2k+1,2k+3} = v_{k+1}^T A v_{k} 
 = - \bigfrac{\delta_{-k-1} \delta_{k+1}}{\beta_{k}}
}}
and
\[
p_{i,2k+1} = p_{2k+1,i} = 
 v_i^T Av_{k} = 0 \tim{for $i \neq -k-1, k-1,k, -k$ and $k+1$,}
\]
\revised{where we note that $\beta_k = v_k^T A^{-1} v_k > 0$ as $A$ 
is positive definite.}
Thus, it follows from from \req{bAvk} and \req{even} that
\eqn{Avk}{\arr{rl}{A v_{k} \ngap  & = 
  p_{2k-1,2k+1} v_{k-1} + p_{2k,2k+1} v_{-k} + p_{2k+1,2k+1} v_{k} \;+ \\
 & \gap\gap 
  p_{2k+2,2k+1} v_{-k-1} + p_{2k+3,2k+1} v_{k+1} \tim{for $k = 1, \ldots, \kmax-1,$} 
}}
where now formally any $p_{2k+1,i} = 0 $ if $i \leq 0$.
Moreover, note that
\[
p_{2k,2k+1} = v_{-k}^T A v_{k} = \delta_{k}
\tim{and}
p_{2k+2,2k+1} = v_{-k-1}^T A v_{k} = \alpha_{k}
\]
because of \eqref{odd}, and
\[
- \frac{\beta_{-k} \alpha_{k-1}}{\beta_{k}}
=
- \frac{\delta_{-k} \delta_{k}}{\beta_{k-1}}
\]
by comparing $p_{2k-1,2k+1}$ with $p_{2k-3,2k-1}$ for subsequent $k$.
%\[
% \delta_{k} = - \bigfrac{\beta_{-k} \alpha_{-k}}{\beta_{k}}
%\tim{and}
%\alpha_{k} = - \bigfrac{\delta_{-k-1} \alpha_{-k-1}}{\beta_{k}}
%\]
Finally, if follows from \eqref{it-zero}, \eqref{Amv} with $k=1$ and the
orthogonality of the $v_i$, that
%\footnote{When $\delta_{-1} = 0$, 
%$p_{1,1} = 1 / \beta_{0}$ and $p_{1,2} = p_{1,3} = 0$, and thus 
%$\alpha_{0}$, $\alpha_{-1}$ and $\delta_{1}$ need not be calculated.}
\eqn{pstart}{p_{1,1} =  v_{0}^T A v_{0} = \frac{1 - \delta_{-1} \alpha_{0}}{\beta_{0}}, \;\;
p_{2,1} = v_{-1}^T A v_{0} = - \frac{\delta_{-1} \alpha_{-1}}{\beta_{0}} \equiv 
 \alpha_{0}, \;\;
p_{3,1} = v_{1}^T A v_{0} = - \frac{\delta_{-1} \delta_{1}}{\beta_{0}},
}
where the alternative definition of $p_{1,2}$ 
follows by symmetry and from \eqref{odd} with $k=0$.

Thus $P_{2k+3}$ is a symmetric, pentadiagonal matrix of the form
%{\small
{\footnotesize
\[
P_{2k+3} =
\smatr{ccccccccc}{p_{1,1} & \alpha_0 & p_{3,1} & 0 & 0 & 0 & 0 & \cdots\\
\alpha_0 & \alpha_{-1} & \delta_{1} & 0 & 0 & 0 & 0 & \cdots \\
p_{3,1} & \delta_{1} & p_{3,3} & \alpha_{1} & p_{5,3} & 0 & 0 & \cdots \\
0 & 0 & \alpha_{1} & \alpha_{-2} & \delta_{2} & 0 & 0 & \cdots \\
0 & 0 & p_{5,3} & \delta_{2} & p_{5,5} & \alpha_{2} & p_{7,5} & \cdots \\
0 & 0 & 0 & 0 & \alpha_{2} & \alpha_{-3} & \delta_{3} & \cdots \\
0 & 0 & 0 & 0 & p_{7,5} & \delta_{3} & p_{7,7} & \cdots \\
\vdots & \vdots & \vdots & \vdots & \vdots & \vdots & \vdots & \ddots
} 
\]
\eqn{P}{ \equiv
\smatr{cccccc:c:cc}{ 
p_{1,1} & \cdots & 0 &  0 & 0 & 0 & 0 & 0 & 0 \\
\vdots & \ddots & \ldots & \ldots & \ldots & 
\ldots & \ldots & \ldots & \ldots \\
0 & \vdots & \ddots & \alpha_{k-2} & p_{2k-1,2k-3} & 0 & 0 & 0 & 0  \\
0 & \vdots & \alpha_{k-2} & \alpha_{-k+1} & \delta_{k-1} & 0 & 0 
& 0 & 0 \\
0 & \vdots & p_{2k-1,2k-3} & \delta_{k-1} & 
p_{2k-1,2k-1} & \alpha_{k-1} & p_{2k-1,2k+1} & 0 & 0 \\
0 & \vdots & 0 & 0 & \alpha_{k-1} & \alpha_{-k} & \delta_{k} & 0 & 0 \\
\cdashline{1-6}
0 & \vdots & 0 & 0 & p_{2k+1,2k-1} & \multicolumn{1}{c}{\delta_{k}} & 
p_{2k+1,2k+1} & \alpha_{k} & p_{2k+1,2k+3} \\
\cdashline{1-7}
0 & \vdots & 0 & 0 & 0 & \multicolumn{1}{c}{0} & 
 \multicolumn{1}{c}{\alpha_{k}} & \alpha_{-k-1} & \delta_{k+1} \\
0 & \vdots & 0 & 0 & 0 & \multicolumn{1}{c}{0} & 
 \multicolumn{1}{c}{p_{2k+3,2k+1}} &  \delta_{k+1} & p_{2k+3,2k+3} \\
},
}
}

\noindent
where the entries $p_{2k+1,2k+1}$ and $p_{2k+1,2k-1}$ (etc) are as given above;
the leading sub-matrices $P_{2k}$ and $P_{2k+1}$ are the entries in the
two ``hashed'' boxes. In other words, we may grow the (lower-triangular parts
of the) submatrices 
$P_{i}$ of $P_{2m+1}$, $i = 1, \ldots, 2m+1$, simply by introducing, for each
successive column, the entries
\eqn{podd}{\arr{rl}{p_{2k+1,2k+1} = &
\bigfrac{1 - \beta_{-k} \delta_{k} - \delta_{-k-1} \alpha_{k}}{\beta_{k}},
\;\;  p_{2k+2,2k+1} = \alpha_{k} \tim{and} \\
p_{2k+3,2k+1} = & - \bigfrac{\delta_{-k-1} \delta_{k+1}}{\beta_{k}} 
\tim{for odd $i =2k+1$}}}
and
\eqn{peven}{p_{2k,2k} = \alpha_{-k}, \;\;  p_{2k+1,2k} = \delta_{k} \tim{and}
p_{2k+2,2k} = 0 \tim{for even $i =2k$,}}
starting with the first column \eqref{pstart} when $i=1$.

Trivially, it also follows from the orthogonality of $V_k$ that 
\eqn{vtb}{V^T_{2m+1} b = \|b\| e_1,} 
where $e_1$ is the first unit vector, and thus $b_k =  \|b\| e_1$ in \req{trs}
for this choice of $V_k$.

We may also combine \req{Amvk} and \req{Avk} to deduce that
%A v_{k-1} = p_{2k-3,2k-1} v_{k-2} + p_{2k-2,2k-1} v_{-k+1} 
% + p_{2k-1,2k-1} v_{k-1} + p_{2k,2k-1} v_{-k} + p_{2k+1,2k-1} v_{k}
%Av_{-k} = p_{2k-1,2k} v_{k-1} + p_{2k,2k} v_{-k} + p_{2k+1,2k} v_{k}
\eqn{AV}{
AV_{2k} = V_{2k} P_{2k} + p_{2k+1,2k-1} v_{k} e_{2k-1}^T + p_{2k+1,2k} v_{k} e_{2k}^T}
%A V_{2k} = V_{2k} P_{2k} + p_{2k+1,2k} v_{k} e_{2k}^T}
and 
%Av_{-k} = p_{2k-1,2k} v_{k-1} + p_{2k,2k} v_{-k} + p_{2k+1,2k} v_{k}
%A v_{k} = p_{2k-1,2k+1} v_{k-1} + p_{2k,2k+1} v_{-k} + p_{2k+1,2k+1} v_{k} +
%  p_{2k+2,2k+1} v_{-k-1} + p_{2k+3,2k+1} v_{k+1}
\eqn{AVp}{AV_{2k+1} = V_{2k+1} P_{2k+1} + 
 ( p_{2k+2,2k+1} v_{-k-1} + p_{2k+3,2k+1} v_{k+1} ) e_{2k+1}^T}
%A V_{2k+1} = V_{2k+1} P_{2k+1} + p_{2k+2,2k+1} v_{-k-1} e_{2k}^T 
% + p_{2k+3,2k+1} v_{k+1} e_{2k+1}^T.}
These identities will shortly be used to develop suitable stopping rules for
our algorithm.

It is perhaps worth mentioning that lines 4--5 and the first part of 8, 
9--11 and the first part of 13 and 14--15, 17 and the first part of 8
in Algorithm~\ref{algorithm-eks} are mathematically equivalent to 
\[\arr{rl}{ \beta_{0} & \ngap \define u^T v_0, \\
  \delta_{-1} & \define \sqrt{\|u\|^2 - \beta_{0}^2},\\
  v_{-1} & \define [u - \beta_{0} v_{0}]/\delta_{-1} \\ \hspace*{1mm} \\
    \alpha_{k-1} & \ngap \define u^T v_{k-1}, \alpha_{-k} \define u^T v_{-k}, \\
    \delta_{k} & \define \sqrt{\|u\|^2 - \alpha_{k-1}^2 - \alpha_{-k}^2}, \\
     v_{k} & \define [u - \alpha_{k-1} v_{k-1} - \alpha_{-k} v_{-k} ]/\delta_{k}
}
\]
and
\[\arr{rl}{
    \beta_{-k} & \ngap \define u^T v_{-k}, \beta_{k} \define u^T v_{k}^T, \\
    \delta_{-k-1} & \define \sqrt{\|u\|^2 - \beta_{-k}^2 - \beta_{k}^2}, \\
       v_{-k-1} & \define [u - \beta_{-k} v_{-k} - \beta_{k} v_{k}] /\delta_{-k-1}.
}
\]
However, numerical experience suggests that this form is unstable, 
in a manner akin to that of Gram-Schmidt orthogonalisation 
of a set of vectors and how such instability may be cured by its modified form 
\cite[\S5.2.7--5.2.8]{GoluvanL96}.

\numsection{Solving the convex trust-region subproblem}
\label{definite-subproblem-solve-section}
\setcounter{algocf}{0}

Having determined an orthogonal basis $V_{\ell}$ of $\calK_{\ell}(A,b)$
using the {\tt EKS} algorithm, we now turn to the trust-region subproblem
\req{trs} 
%(with $\ell = 2k+1$) 
within this subspace.
%\footnote{We might 
%also consider the basis $V_{2k}$ for the $\calK_{2k}(A,b)$, and
%the corresponding subproblem \req{trs} (with $\ell = 2k$), but for
%brevity here we focus on $V_{2k+1}$. The $2k$ case is identical.}.
As we have just seen, the Hessian $P_{\ell} = V_{\ell}^T A V_{\ell}$ 
is symmetric and sparse (pentadiagonal),
while the first-order term $b_{\ell}$ is a scalar multiple
$\|b\|$ of the $\ell$-dimensional first unit vector. Significantly, we
are expecting $\ell$ to be small, and thus the subproblem is low dimensional.

Thus we might apply off-the-shelf trust-region software such as 
{\tt GQT} \cite{MoreSore83:sisc} or
{\tt TRS} \cite{GoulRobiThor10:mpc} from 
{\sf GALAHAD} \cite{GoulOrbaToin03:toms}.
However, the alternative we prefer is to proceed as in \req{spectral}, that is
to compute the spectral factorization
\eqn{spectralk}{P_{\ell} = U_{\ell} \Lambda_{\ell} U_{\ell}^T} 
involving matrices of eigenvectors $U_{\ell}$ and (diagonal) eigenvalues 
$\Lambda_{\ell}$, and to find $y_{\ell} = U_{\ell} z_{\ell}$, where
\eqn{trtk}{z_{\ell} = \argminin{z\in \tinyRe^{\ell}} 
 \half z^T \Lambda_{\ell} z - d_{\ell}^T z \tim{subject to}  \|z\| \leq \Delta}
and $d_{\ell} = U_{\ell}^T b_{\ell}$.
Significantly \req{spectralk} is inexpensive since $\ell$ is small, and
the Hessian of \req{trtk} is diagonal. This then enables very fast and
accurate solution using high-order, factorization-free iterative schemes, 
such as {\tt TRS\_solve\_diagonal} \cite{GoulRobiThor10:mpc} 
from {\sf GALAHAD}. 

Having computed the solution $y_{\ell}$ to \req{trs} as above 
(or by any other means), 
together with the corresponding shift $\sigma_{\ell}$ for which both
\eqn{yopt}{( P_{\ell} + \sigma_{\ell}I)y_{\ell} =  \|b\| e_1}
and $\|x_{\ell}\| = \Delta$ hold, according to \req{errory}, 
it remains to see how big is the residual 
\disp{r_{\ell} \eqdef ( A + \sigma_{\ell}I)x_{\ell} - b,}
where $x_{\ell} = V_{\ell} y_{\ell}$. When $\ell = 2k+1$, we find that
\[
\arr{rl}{A x_{2k+1} \ngap
 & = A V_{2k+1} y_{2k+1} \\
 & = V_{2k+1} P_{2k+1} y_{2k+1} +
 ( p_{2k+2,2k+1} v_{-k-1} + p_{2k+3,2k+1} v_{k+1} ) e_{2k+1}^Ty_{2k+1} \\
 & = V_{2k+1} (  \|b\|e_1 - \sigma_{2k+1} y_{2k+1}) +
 ( p_{2k+2,2k+1} v_{-k-1} + p_{2k+3,2k+1} v_{k+1} ) e_{2k+1}^Ty_{2k+1} \\
 & = b - \sigma_{2k+1} x_{2k+1} +
 ( p_{2k+2,2k+1} v_{-k-1} + p_{2k+3,2k+1} v_{k+1} ) e_{2k+1}^Ty_{2k+1} \\
}
\]
%     p_{2k+2,2k+1} v_{-k-1} e_{2k}^T y_{2k+1} + p_{2k+3,2k+1} v_{k+1} e_{2k+1}^T y_{2k+1} \\
% & = V_{2k+1} (  \|b\|e_1 - \sigma_{2k+1} y_{2k+1}) +
%     p_{2k+2,2k+1} v_{-k-1} e_{2k}^T y_{2k+1} + p_{2k+3,2k+1} v_{k+1} e_{2k+1}^T y_{2k+1} \\
% & = b - \sigma_{2k+1} x_{2k+1} +
%     p_{2k+2,2k+1} v_{-k-1} e_{2k}^T y_{2k+1} + p_{2k+3,2k+1} v_{k+1} e_{2k+1}^T y_{2k+1} \\
%}
%\]
using \req{AVp}, \req{yopt}, the orthogonality of the columns of $V_{2k+1}$, 
and the relationship $ V_{2k+1} \|b\| e_1 = b$.
Thus 
%\[r_{2k+1} = p_{2k+2,2k+1} v_{-k-1} e_{2k+1}^T y_{2k+1} + 
%            p_{2k+3,2k+1} v_{k+1} e_{2k+1}^T y_{2k+1}
\[r_{2k+1} = ( p_{2k+2,2k+1} v_{-k-1} + p_{2k+3,2k+1} v_{k+1} ) e_{2k+1}^Ty_{2k+1}
\]
and hence
\eqn{rkp}{\|r_{2k+1}\| = |y_{2k+1,2k+1}| \sqrt{p^2_{2k+2,2k+1} + p^2_{2k+3,2k+1}},}
where $y_{2k+1,2k}$ and $y_{2k+1,2k+1}$ are the $2k$-th and $2k+1$-st
components of $y_{2k+1}$, using the orthogonality of $v_{-k-1}$ and $v_{k+1}$.
Thus the norm of the residual may be computed at almost no extra cost
using data that is already available.
Identical reasoning, when $\ell = 2k$ but now using \req{AV}, shows that
%\[r_{2k} = p_{2k+1,2k} v_{-k-1} e_{2k-1}^T y_{2k} + 
%          p_{2k+2,2k} v_{k+1} e_{2k}^T y_{2k}
\[r_{2k} = ( p_{2k+1,2k-1} e_{2k-1}^T y_{2k} + p_{2k+1,2k} e_{2k}^T y_{2k} ) v_{k}\]
and hence in this case
\eqn{rk}{\|r_{2k}\| = | p_{2k+1,2k-1} y_{2k,2k-1} + p_{2k+1,2k} y_{2k,2k} |.}
Note that in the exceptional case where $\delta_{\ell} = 0$, the corresponding
$\|r_{\ell}\| = 0$, and thus termination will always occur at this step;
\revised{the {\tt EKS} algorithm has generated the complete extended Krylov 
subspace, and the solution recovered from \req{trs} will be that of
\req{tr} (in everything but the probability-zero hard case).}

The {\tt TREK} algorithm embeds solutions of \req{trs} after 
every computation of $v_{-k}$ and $v_k$ in {\revised{Algorithm~\req{algorithm-eks}}}; at these
stages the matrices $V_{2k}$ and $V_{2k+1}$, respectively, are complete,
and we have the data required to solve the appropriate subproblems \req{trs}.
Armed with such solutions, we terminate as soon as $\|r_{2k}\|$ from 
\req{rk} or $\|r_{2k+1}\|$ from \req{rkp} is deemed small enough.
One further important point is that after computing 
$u \define A^{-1} v_{0} = A^{-1} b / \|b\|$ before we start the main loop in the 
{\tt EKS} algorithm, it is trivial to check for the interior
case \req{trivial}. Thus if we check and find 
$\|u\| \|b\| \leq \Delta$ in {\tt TREK},
the solution is interior, and we exit\footnote{Actually, we compute 
$x = A^{-1} b$, exit with $x_* = x$ if $\|x\| \leq \Delta$, and otherwise
continue with $u = x /\|b\|$.} with $x_* = \|b\| u$. For completeness
we state the full {\tt TREK} algorithm (Algorithm~\ref{algorithm-trek})
\vpageref{algorithm-trek}.

\begin{algorithm2e}[!ht]
\caption{The {\tt TREK} algorithm to solve the trust-region subproblem
\req{tr} \label{algorithm-trek}}
\KwInput{symmetric $A \in \Re^{n\times n}$, $b \in \Re^n$, $\Delta > 0$,
 a stopping tolerance $\epsilon > 0$  and an iteration bound $m\geq 1$.}
\KwOutput{$x_* = \approxargmininx{} \half x^T A x - b^T x 
 \tim{subject to} \|x\| \leq \Delta.$}
  \setstretch{0.8}
  \vspace{2mm}
  {\bf Algorithm:}\\
  $x \define A^{-1} b$\;
  \If{$\|x\| \leq \Delta$}{
  {\bf exit} with the interior solution $x_* = x$,
    shift $\sigma_* = 0$ and $k_{\mbox{\scriptsize stop}} = 0$\;
  }
  $\delta_0 \define \|b\|$, 
  $v_{0} \define b/\delta_0$, 
  $u \define x/\delta_0$\;
  $\beta_{0} \define u^T v_0$, 
  $u \redef u - \beta_{0} v_{0}$\;
  $\delta_{-1} \define \|u\|$\;
  \For{$k = 1, 2, \ldots , m$}{
    \uIf{$\delta_{-k} > 0$}{
     $v_{-k} \define u/\delta_{-k}$,
     $u \define Av_{-k}$\;
     $\alpha_{k-1} \define u^T v_{k-1}$,
     $u \redef u - \alpha_{k-1} v_{k-1}$\;
     $\alpha_{-k} \define u^T v_{-k}$,
     $u \redef u - \alpha_{-k} v_{-k}$\;
     $\delta_{k} \define \|u\|$\;
    }
    \uIf{$k>1$}{update $P_{2k-1}$ from $P_{2k-2}$ using \req{podd} and form
     $b_{2k-1} = (b^T_{2k-2},0)^T$\;}
    \Else{initialize $P_1$ from \req{pstart} and set $b_1 = \delta_0$\;}
%   update $P_{2k-1}$ from $P_{2k-2}$ from \req{P} and form
%    $b_{2k-1} = (b^T_{2k-2},0)^T$\;
    $y_{2k-1} \define \argmin \half y^TP_{2k-1} y - b_{2k-1}^Ty \; \mbox{s.t.} \;
      \|y\| \leq \Delta$ and its optimal shift $\sigma_{2k-1}$\;
    compute $\|r_{2k-1}\|$ from \req{rkp}\;
    \If{$\|r_{2k-1}\| \leq \epsilon$}{
    {\bf exit} with the interior solution $x_* = V_{2k-1} y_{2k-1}$, shift
      $\sigma_* = \sigma_{2k-1}$ and $k_{\mbox{\scriptsize stop}} = k$\;}
    \uIf{$\delta_{k} > 0$}{
     $v_{k} \define u/\delta_{k}$, 
     $u \define A^{-1} v_{k}$\;
     $\beta_{-k} \define u^T v_{-k}$,
     $u \redef u - \beta_{-k} v_{-k}$\;
     $\beta_{k} \define u^T v_{k}$,
     $u \redef u - \beta_{k} v_{k}$\;
    }
    \If{$k<m$}{\uIf{$\delta_{k} > 0$}{$\delta_{-k-1} \define \|u\|$\;}
%    update $P_{2k+1}$ from $P_{2k}$ from \req{P} and form
%    $b_{2k+1} = (b^T_{2k},0)^T$\;
%    $y_{2k+1} \define \argmin \half y^TP_{2k+1} y - b_{2k+1}^Ty \; \mbox{s.t.} \;
%      \|y\| \leq \Delta$ and its optimal shift $\sigma_{2k+1}$\;
%    compute $\|r_{2k+1}\|$ from \req{rkp}\;
%    \If{$\|r_{2k+1}\| \leq \epsilon$}{
%    {\bf exit} with the interior solution $x_* = V_{2k+1} y_{2k+1}$, shift
%      $\sigma_* = \sigma_{2k+1}$ and $k_{\mbox{\scriptsize stop}} = k+1$\;
%    }
     update $P_{2k}$ from $P_{2k-1}$ using \req{peven} and form
      $b_{2k} = (b^T_{2k-1},0)^T$\;
%     \uIf{$k>1$}{update $P_{2k}$ from $P_{2k-1}$ from \req{P} and form
%      $b_{2k} = (b^T_{2k-1},0)^T$\;}
%     \Else{initialize $P_{2}$ from \req{P} and set 
%      $b_{2} = \delta_0 e_1$\;}
     $y_{2k} \define \argmin \half y^TP_{2k} y - b_{2k}^Ty \; \mbox{s.t.} \;
       \|y\| \leq \Delta$ together with its optimal shift $\sigma_{2k}$\;
     compute $\|r_{2k}\|$ from \req{rk}\;
     \If{$\|r_{2k}\| \leq \epsilon$}{
     {\bf exit} with the interior solution $x_* = V_{2k} y_{2k}$, shift
       $\sigma_* = \sigma_{2k}$ and $k_{\mbox{\scriptsize stop}} = k$\;
     }
  }
 }
\end{algorithm2e}
\vspace*{-2mm}

%\begin{itemize}
%\item estimate errors from {\tt EKS} subspace approximation 
%$x_k(\lambda) = V_k y_k(\lambda)$ with
%$(P_k + \lambda I ) y_k(\lambda) = \|b\| e_1$ to 
%$(A + \lambda I) x(\lambda) = b$
%results when%using the error estimate $| \gamma_k(\lambda)|$ where
%$(A + \lambda I) x(\lambda) - b = v_k \gamma_k(\lambda)$.
%\item for a series of ordered $\lambda_i$ (say, smallest to largest), 
%use the above and truncate search for the first $i$ for which 
%$| \gamma_k(\lambda)| > \epsilon$.
%\item order search using one of the following: 
%(i) from smallest to largest or vice versa 
%(ii) pick initial $i$ to give the mean of the $\lambda_i$,
%and search ``outwards'' from this
%(iii) pick initial $i$ to give the median of the $\lambda_i$,
%and search ``outwards'' from this
%\item on subsequent searches start with the value of $i$ at which
%failure occurred.
%\item can monitor every iteration, every other or at some other interval.
%\end{itemize}

\subsection{Resolves}
\label{resolves-section}

The need to resolve a problem with the same $A$ and $b$, but a smaller
$\Delta$ is a common occurrence in trust-region methods, since such a
mechanism is used to guarantee convergence when the current solution
has proved unsatisfactory \cite{ConnGoulToin00}. This is easy to
accommodate within Algorithm~\ref{algorithm-trek}; the data generated
by the previous call to the algorithm is retained,
and the ``for loop'' re-entered with $k$ taking the value 
$k_{\mbox{\scriptsize stop}}$ that it last
had. The motivation is simply that although the 
\revised{subspace $\calK_{2k}(A,b)$ that we have used to approximate 
$\calS(A,b)$}
% \sout{segment of the manifold}} 
%we have observed}
\revised{might now  need} to expand slightly,
it may well be that there is little actual difference; the data provided
in Appendix~\ref{sec:detailed numerical experiments} indicates that in 
only \revised{six of ninety three} cases examined 
were any further extended Krylov iterations necessary, and thus that
resolves are generally very inexpensive. 
Examining the bound \req{rate} given in Proposition~\ref{prop:convergence_rate} in 
Appendix~\ref{sec:convergence analysis}, and looking at Figure~\ref{fig:E},
we see that the convergence rate of the extended Krylov subspace
improves for $\sigma$ beyond a certain point, which is consistent with this
observation.
Although other methods can also take advantage of resolves, all that 
we are aware of still require significant extra computation.

\subsection{Other elliptical norms}
\label{other-elliptical-norms-section}

Although the trust-region is often defined in terms of the Euclidean norm,
it can be that this promotes some variables at the expense of others. To
cope with this, it is sometimes more effective to pose the problem as
\eqn{trp}{
 \mininx{} \half x^T A x - b^T x \tim{subject to} \|\PP x\| \leq \Delta,}
where $\PP$ is an easily-invertible matrix---a diagonal $\PP$
is often a simple but effective choice. On writing 
\eqn{xp}{\xp = \PP x,}
the minimizer of \req{trp} is formally $x_* = \PP^{-1} \xpstar$, where
\eqn{trps}{\xpstar = 
 \argmininx{} \half \xp{}^T \Ap \xp - \bp{}^T \xp 
 \tim{subject to} \|\xp\| \leq \Delta}
and $\Ap = \PP^{-T} A \PP^{-1}$ and $\bp = \PP^{-T} b$. 
Since the {\tt TREK} algorithm
described above merely relies on products with $\Ap$ and its inverse
$\Ap{}^{-1} = \PP A^{-1} \PP^T$, once again we merely need a factorization of $A$
(to find $A^{-1}v$ for given $v$) so long as means of inverting $\PP$ and 
its transpose are available. Commonly a positive-definite scaling 
matrix $S$ that is (in some sense) close to $A$ is known, and in this
case a decomposition, for example a sparse Cholesky factorization, 
$S = \PP^T \PP$ gives the appropriate $\PP$. We provide details of the necessary
changes to devise suitable algorithms to solve the resulting problem 
\eqn{trss}{
 \mininx{} \half x^T A x - b^T x \tim{subject to} \|x\|_S \leq \Delta,}
where $\|x\|_S \eqdef \sqrt{x^T S x}$, in 
Appendix~\ref{sec:modified norm subproblem solver}. We illustrate one such
algorithm here as Algorithm~\ref{algorithm-treks2} \vpageref{algorithm-treks2}.
Note that here the generated vectors $v_i$ now form an $S$-orthonormal basis
for the extended Krylov subspace.

\begin{algorithm2e}[!ht]
\caption{The {\tt TREK} algorithm to solve the trust-region subproblem
\req{trss} \label{algorithm-treks2}}
\KwInput{symmetric $A \in \Re^{n\times n}$, $b \in \Re^n$, 
 symmetric, positive-definite $S \in \Re^{n\times n}$, $\Delta > 0$,
 a stopping tolerance $\epsilon > 0$  and an iteration bound $m\geq 1$.}
\KwOutput{$x_* = \approxargmininx{} \half x^T A x - b^T x 
 \tim{subject to} \|x\|_S \leq \Delta.$}
  \setstretch{0.8}
  \vspace{2mm}
% \small
  {\bf Algorithm:}\\
  $x \define A^{-1} b$, 
  $u \define S x$\;
  \If{$\sqrt{u^T x} \leq \Delta$}{
  {\bf exit} with the interior solution $x_* = x$,
    shift $\sigma_* = 0$ and $k_{\mbox{\scriptsize stop}} = 0$\;
  }
  $w \define S^{-1} b$, 
  $\delta_0 \define \sqrt{b^T w}$\;
  $v_0 \define w / \delta_0$, 
  $q_0 \define b/\delta_0$, 
  $w \define u / \delta_0$,
  $u \define x / \delta_0$\;
  $\beta_0 \define u^T q_0$,
  $u \redef u - \beta_0 v_0$\;
  $w = S u$, 
  $\delta_{-1} \define \sqrt{u^T w}$\;
  \For{$k = 1, 2, \ldots , m$}{
    \uIf{$\delta_{-k} > 0$}{
     $v_{-k} \define u/\delta_{-k}$,
     $q_{-k} \define w/\delta_{-k}$\;
     $w \define A v_{-k}$, 
     $u \define S^{-1} w$\;
     $\alpha_{-k} \define w^T v_{-k}$,
     $u \redef u - \alpha_{-k} v_{-k}$,
     $w \redef w - \alpha_{-k} q_{-k}$\;
     $\alpha_{k-1} \define w^T v_{k-1}$,
     $u \redef u - \alpha_{k-1} v_{k-1}$,
     $w \redef w - \alpha_{k-1} q_{k-1}$\;
     $\delta_{k} \define \sqrt{u^T w}$\;
    }
    \uIf{$k>1$}{update $P_{2k-1}$ from $P_{2k-2}$ using \req{podd} and form
     $b_{2k-1} = (b^T_{2k-2},0)^T$\;}
    \Else{initialize $P_1$ from \req{pstart} and set $b_1 = \delta_0$\;}
    $y_{2k-1} \define \argmin \half y^TP_{2k-1} y - b_{2k-1}^Ty \; \mbox{s.t.} \;
      \|y\| \leq \Delta$ and its optimal shift $\sigma_{2k-1}$\;
    compute $\|r_{2k-1}\|$ from \req{rkp}\;
    \If{$\|r_{2k-1}\| \leq \epsilon$}{
    {\bf exit} with the interior solution $x_* = V_{2k-1} y_{2k-1}$, shift
      $\sigma_* = \sigma_{2k-1}$ and $k_{\mbox{\scriptsize stop}} = k$\;}
    \uIf{$\delta_{k} > 0$}{
     $v_{k} \define u/\delta_{k}$,
     $q_{k} \define w/\delta_{k}$\;
     $u \define  A^{-1} q_k$,
     $w \define  S u$\;
     $\beta_{-k} \define w^T v_{-k}$,
     $u \redef u - \beta_{-k} v_{-k}$, 
     $w \redef w - \beta_{-k} q_{-k}$\;
     $\beta_{k} \define w^T v_{k}$,
     $u \redef u - \beta_{k} v_{k}$,  
     $w \redef w - \beta_{k} q_{k}$\;
    }
    \If{$k<m$}{\uIf{$\delta_{k} > 0$}{$\delta_{-k-1} \define \sqrt{u^T w}$\;}
     update $P_{2k}$ from $P_{2k-1}$ using \req{peven} and form
      $b_{2k} = (b^T_{2k-1},0)^T$\;
     $y_{2k} \define \argmin \half y^TP_{2k} y - b_{2k}^Ty \; \mbox{s.t.} \;
       \|y\| \leq \Delta$ together with its optimal shift $\sigma_{2k}$\;
     compute $\|r_{2k}\|$ from \req{rk}\;
     \If{$\|r_{2k}\| \leq \epsilon$}{
     {\bf exit} with the interior solution $x_* = V_{2k} y_{2k}$, shift
       $\sigma_* = \sigma_{2k}$ and $k_{\mbox{\scriptsize stop}} = k$\;
     }
  }
 }
\end{algorithm2e}
%\noindent

It is also worth noting that the
extended Krylov subspace \req{eks} generated by  $\Ap$ and $\bp$,
\eqn{peks}{\calK_{2m+1}(\Ap,\bp) = 
 \text{span}\{\bp,\Ap^{-1}\bp,\Ap \bp,\ldots,\Ap^{-m}\bp,\Ap^m\bp\},}
is
\disp{\PP^{-T}\text{span}\{b,(SA^{-1})b,(AS^{-1}) b,\ldots,\revised{(SA^{-1})^m}b,(AS^{-1})^m b\},}
and it is then possible to derive a suitable variant of 
Algorithm~\ref{algorithm-eks} using $S$ and its inverse rather 
than $\PP$ and $\PP^{-1}$.

\clearpage
\numsection{Solving the general trust-region subproblem}
\label{general-subproblem-solve-section}

For non-convex problems---indeed if $A$ is not positive definite---it may 
be awkward (or even impossible if $A$ is singular) 
to use an extended Krylov subspace involving 
products with $A^{-1}$. Fortunately, it is straightforward to deal with this
possibility using the subspace generated by a
diagonally-shifted matrix $\As = A + \sigmas I$, 
for some suitable $\sigmas > 0$, and its inverse. 

To see why, first note that the interior-solution case cannot occur if
$A$ is indefinite, and thus the solution we seek must lie on the trust-region
boundary. That is, we must satisfy \req{ls}. We may rewrite this as
\eqn{lss}{ ( \As + \theta I ) \xt = b \tim{with} \|\xt\| = \Delta,}
where $\xs \equiv \xt$ and
\disp{\theta = \sigma - \sigmas
\geq \sigmamin - \sigmas = \max(- \sigmas,- \lambda_1(A)- \sigmas).}
Thus we seek the unique root of \req{lss} but now over the interval
$\theta \in [\sigmamin - \sigmas,\infty)$.
Given $\sigmas$, it is trivial to adapt any bracketing root-finding 
procedure to account for this extended interval.

We have in mind that the shift should \revised{be chosen} 
%\revised{\sout{be}} 
so that $\As$ is positive definite. The smallest shift for which 
$\As$ is guaranteed to be positive semi-definite is $\sigmamin$ from
\req{sigmamin} so any value slightly larger than this is possible.
But a far cheaper option that does not rely on knowing the
spectrum of $A$, once $A$ has been found not to be positive definite,
is to use the value 
\disp{\sigma\sub{G} = - \min_{1\leq i \leq n} \left[ \revised{[A]_{i,i}} 
 - \sum_{1\leq j \neq i \leq n} | \revised{[A]_{i,j}}| \right] + \epsilon\sub{G}}
from the Gershgorin theorem, requiring a single pass through the values
of $A$; we have found that this with the heuristic value 
$\epsilon\sub{G} = \sqrt{\epsilon\sub{M}} \max_{1\leq i,j \leq n} 
\revised{|[A]_{i,j}}|$
where $\epsilon\sub{M}$ is the machine precision is effective.

As before, we generate the extended Krylov subspace $\calK_{\ell}(\As,b)$ using
Algorithm~\ref{algorithm-eks}. If the $\sigma$ we seek lies to the right
of  $\sigmas$ (i.e., $\sigma \geq \sigmas$ and hence $\theta \geq 0$), the
Beckermann-Townsend \cite[\S4.2]{BeckTown17:simaa} analysis
outlined in Section~\ref{low-rank-section} continues to hold, and we can
expect that $\xt$ will lie in an approximately low-dimensional subspace that
will be uncovered by Algorithm~\ref{algorithm-eks}. 

It is also worth noting that although we have now generated
$P_{\ell} = V_{\ell}^T \As V_{\ell}$, we equally have
\[ V_{\ell}^T A V_{\ell} =  V_{\ell}^T \As V_{\ell} -  \sigmas V_{\ell}^T V_{\ell} 
 = P_{\ell} - \sigmas I\]
since the columns of {\revised{$V_{\ell}$}} are orthonormal. Thus subspace
trust-region subproblems of the form \req{trs} can be formulated as
\eqn{trsns}{y_{\ell} =  \revised{\argminin{y \in \tinyRe^{\ell}}} 
 \half y^T ( P_{\ell} - \sigmas I ) y - b_{\ell}^T y  \tim{subject to} 
\|y\| \leq \Delta.}
More particularly, if we rely on the spectral factorization \req{spectral},
the diagonal subproblem \req{trtk} becomes
\[z_{\ell} = \argminin{z\in \tinyRe^{\ell}} 
 \half z^T ( \Lambda_{\ell} - \sigmas I ) z - d_{\ell}^T z 
\tim{subject to} \|z\| \leq \Delta.
\]

Finally, we return briefly to the thorny issue of the so-called hard case, 
that is when $b$ is orthogonal to the subspace $\calE_1$ of all the eigenvectors 
for the leftmost eigenvalue $\lambda_1(A) < 0$, and the radius is too large. 
In this case, the solution we seek is of the form\footnote{Note that
$x_+$ exists since $b$ lies in the range of $A - \lambda_1(A)I$
as $b$ is orthogonal to $\calE_1$.}
\[x_+ + \theta y, \tim{where}
x_+ = \ngap\ngap\lim_{\sigma \stackrel{+}{\rightarrow} - \lambda_1(A)} \ngap\ngap 
x(\sigma), \;\; y \in \calE_1, \tim{and} \|x_+ + \theta y\| = \Delta,
\]
when $ \|x_+\| < \Delta$.
No Krylov method based on $A$ and $b$ (extended or otherwise) will see the 
eigenspace $\calE_1$ in the hard case in exact arithmetic, although rounding 
errors can gradually introduce it. Since we have never observed the hard case 
in practice for anything other than contrived examples, our only precaution is
to add a tiny (pseudo-random) perturbation to $b$ if requested. 
We certainly do this when $b = 0$.

\subsection{Other elliptical norms}

For the general problem defined in terms of an elliptical norm
$\|x\|_S = \sqrt{x^T S x}$ for a given symmetric positive definite 
$S = \PP^T \PP$, a similar approach is possible. Transforming via \req{xp} as
in Subsection~\ref{other-elliptical-norms-section}, the (non-interior) 
optimality condition \revised{\req{ls}} becomes
\eqn{lsp}{ ( \Ap + \sigma I ) \xps = \bp \tim{with} \|\xps\| = \Delta}
for some
\disp{\sigma \geq \sigmapmin \eqdef \max(0,- \lambda_1(\Ap)),}
or equivalently
\disp{ ( A + \sigma S ) \xs = b \tim{with} \|\xs\|_S = \Delta;}
note that $\lambda_1(\Ap)$ is equivalently the leftmost eigenvalue
$\lambda_1(A,S)$ of the symmetric matrix pencil \cite[\S7,7]{GoluvanL96}.
$A - \lambda S$. The equivalent version of \req{lss} is then
\eqn{lsps}{ ( \As + \theta S ) \xt = b \tim{with} \|\xt\|_S = \Delta,}
where 
\disp{\theta = \sigma - \sigmas
\geq \sigmapmin - \sigmas = \max(- \sigmas,- \lambda_1(A,S)- \sigmas).}
Although it is less obvious how to find a lower bound on
$\lambda_1(A,S)$ in general, it is possible to use a variant of the 
inexpensive Gershgorin bound in the commonly-occurring case where $S$ is 
strictly diagonally dominant (see \cite[\S3.3.2]{GoulRobiThor10:mpc}
for details). Our software packages make the strict diagonal dominance
of $S$ a prerequisite.
%\footnote{However, for the time being, we have restricted 
%$P$ to be diagonal matrices, as more general $P$ requires additional external 
%software---in particular access to products with $S$ or $P$---that has yet 
%to be generally implemented.}

\numsection{Numerical experiments}
\label{experiments-section}

We now provide some evidence that our new approach is a useful alternative
to common existing methods. We consider the set of \revised{93} 
unconstrained test
problems in the current release (2025-09-09) of the \cutest\ optimization 
test examples \cite{GoulOrbaToin15:coap}
that have 1000 or more variables. For these, we generate 
the gradient and Hessian at the provided initial point, and use these for $-b$
and $A$ respectively. 

We compare our new method (Algorithm~\ref{algorithm-trek} implemented 
as {\tt TREK}), using the initial (Gershgorin-based) shifting strategy 
outlined in Section~\ref{general-subproblem-solve-section},
with the multi-factorization 
method \cite{GoulRobiThor10:mpc} available as {\tt TRS} 
within \galahad\ \cite{GoulOrbaToin03:toms} 
and the iterative method {\tt GLTR} \cite{GoulLuciRomaToin99:siopt}
from the same library; {\tt GLTR} only requires Hessian-vector
products which can be an asset in some cases. Since all three methods 
can take advantage of information gathered when solving a problem with one
radius and are now faced with another with identical data but a smaller radius,
we consider each problem with \revised{two initial radii, $\Delta = 10$ and 
$\Delta = 1$, and then reduce the radii by a factor of $10$ until
$\Delta = 0.1$ is reached.} {\tt TREK}
uses an iteration bound \revised{$m = 300$}, but almost always\footnote{In 
\revised{six} instances, $m$ was too small, but increasing the value to $400$ 
cured the problem.} 
terminates with far-fewer $k$. Default, and roughly equivalent
stopping rules were applied; {\tt TREK} terminates as soon as
\req{rkp} or \req{rk} is smaller than $10^{-10}$.

We perform our experiments on four cores of a PC with sixteen
Intel Core i9-9900 CPU 3.10GHz processors, and thirty two Gbytes of memory.
The codes are all from the \galahad\ library compiled with 
{\tt -Ofast} optimization
by gfortran, and use is made of OpenMP-tuned BLAS and LAPACK. Factorization is
performed, when needed, by the {\tt HSL\_MA57} sparse symmetric linear solver
\cite{DuffReid96b:toms}, enabled via \galahad's {\tt SLS} package, and
relies on the tuned BLAS for good performance. On some examples, particularly
those whose Hessians are banded, further experiments not listed here 
revealed that LAPACK's {\tt PBTR} band solver is a good alternative.

\begin{figure}[htbp]
\begin{center}
\vspace*{1.5mm}
\includegraphics[height=6cm]{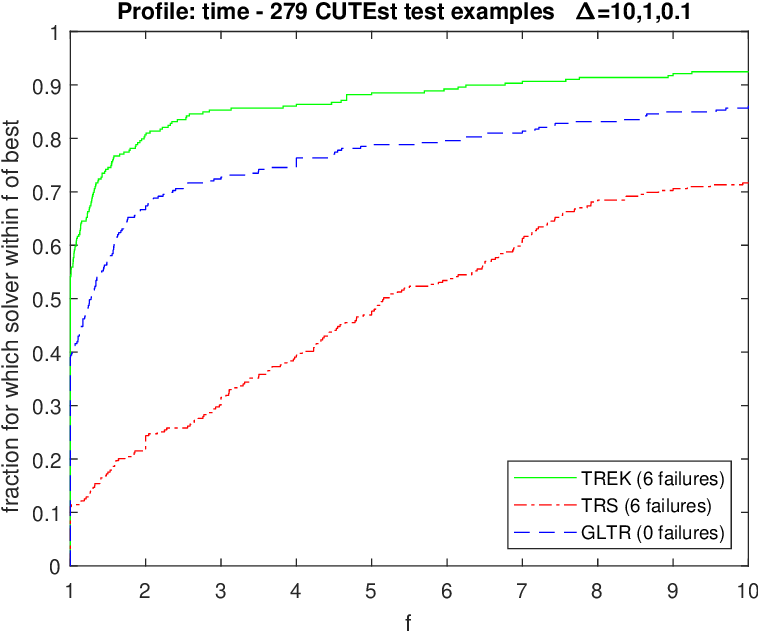} 
\includegraphics[height=6cm]{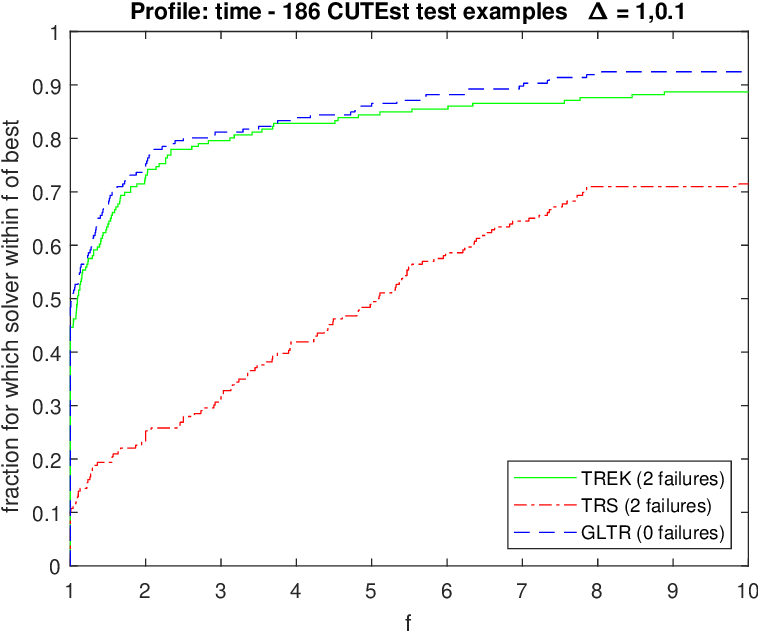} 
% produced by profile__103_time.m and profile__12_time.m
\caption{\label{figure-times} Performance profile comparing three
trust-region subproblem solvers on larger \cutest\ examples.
\revised{Here, for each factor $f$ on the horizontal axis, the vertical axis 
shows the fraction of the number of examples that were solved within a ratio 
$f$ of the best result over all three solvers on that example, 
measured according to CPU time taken. The figure on the left is for runs
started with $\Delta = 10$, 
and the solutions are then updated at $\Delta = 1$ and $0.1$. That
on the right is from $\Delta = 1$ and updated at $\Delta = 0.1$.
}}
\end{center}
\end{figure}

We summarize our findings in Figure~\ref{figure-times} 
using the performance profile \cite{dolamore02:mp}
in which relative times are ranked.
The profile is built using the data given in Appendix~\ref{sec:detailed numerical experiments}.

We did not expect any of the methods tested to be an overall winner,
and indeed this is the case. Clearly, the factorization-free method
{\tt GLTR} is frequently very effective, particularly when factorization
is a significant overhead \revised{and the initial radius is small}. 
But we were encouraged to see that {\tt TREK}
performs well in other cases, and the requirement of a single factorization
often proves decisive when compared to the multi-factorization approach
used by {\tt TRS}.\footnote{{\tt TRS} also failed in \revised{six} cases 
because the shifted matrix $A + \sigma I$ became too close to singularity to
progress.} \revised{Significantly, and as predicted in the 
motivating justification in Section~\ref{manifoild-approximation-section},
the left-hand figure in Figure~\ref{figure-times} seems to indicate that
including $A^{-k} b$ terms in the approximating subspace
is advantageous when $\Delta$ is large 
(i.e., $\sigma$ is closer to $\sigmamin$).}

\numsection{Regularization subproblems}
\label{regularization-section}

Another important subproblem employed by methods for nonlinear optimization
uses norm-regularization 
\cite{CartGoulToin11:mp,CartGoulToin22,NestPoly06:mp}. 
Again, given an $n$ by $n$  symmetric (Hessian) matrix $A$, 
an $n$-vector $b$ and the Euclidean norm $\|\cdot\|$, but now
a regularization weight $\rho > 0$ and a power $r \geq 2$,
the norm-regularization subproblem is to find
\eqn{rn}{x_* = 
 \argmininx{} \half x^T A x - b^T x + \frac{\rho}{r} \|x\|^r.}
Once again, the solution $\xs$ satisfies the first-order optimality 
condition \req{fon}, but now
\eqn{sigmarn}{\sigma = \rho \|\xs\|^{r-2}}
for some $\sigma \geq \sigmamin$
\revised{\cite[Thm.8.2.8]{CartGoulToin22}.} 
Thus the \revised{set}
%\sout{manifold}
of solutions  \revised{$\{x(\sigma)\}$} is identical to that for the 
trust-region subproblem, but the one of interest satisfies \req{sigmarn} 
instead of \req{ls}; there is no interior-solution possibility in this case.

Since this is the case, it is trivial to recast the {\tt TREK} algorithm
for \req{rn}, the only differences are that the steps that calculate
$y_{2k}$ and $y_{2k+1}$ and their shifts find instead
\eqn{nk}{y_{2k} \define \argmin \half y^TP_{2k} y - b_{2k}^Ty + 
 \frac{\rho}{r} \|y\|^r \tim{and its optimal shift $\sigma_{2k}$}}
and
\eqn{nkp}{y_{2k+1} \define \argmin \half y^TP_{2k+1} y - b_{2k+1}^Ty + 
 \frac{\rho}{r} \|y\|^r \tim{and its optimal shift $\sigma_{2k+1},$}}
and that the check for an interior solution before the main
iteration loop is omitted. We refer to the resulting algorithm 
(and the {\sf GALAHAD} package that implements it) as
{\tt NREK}. The subproblems \req{nk} and \req{nkp} are solved just
as those for the trust-region case by finding the spectral
decompositions of $P_{2k}$ and $P_{2k+1}$, and transforming so that the 
resulting (small) subproblem is diagonal. Thereafter the high-order, 
factorization-free diagonal regularization subproblem subroutine 
{\tt RQS\_solve\_diagonal} \cite{GoulRobiThor10:mpc} from {\sf GALAHAD} 
is employed. Other elliptical norms may be handled just as described in
Subsection~\ref{other-elliptical-norms-section}. 

\numsection{Conclusions}
\label{conclusions-section}

We have provided an alternative to currently popular methods for
solving quadratic trust-region and regularization problems. This
lies somewhere between multiple factorization methods that were
first popularised by Mor\'{e} and Sorensen \cite{MoreSore83:sisc}
and CG/Lanczos-based factorization-free methods such as {\tt GLTR}
\cite{GoulLuciRomaToin99:siopt}. The novel features are firstly that
it is shown that the \revised{set} 
%\sout{manifold}
of all possible solutions lies in an (approximately) very 
\revised{low-dimensional subspace},
%\sout{manifold}
and secondly
that by extending the Krylov subspace to include products with $A^{-1}$ as well
as $A$, this \revised{subspace}
%\sout{manifold}
may be computed very efficiently. This then adds
a further tool to the arsenal of solvers that trust-region and 
cubic-regularization methods for general optimization methods rely on.
Unsurprisingly, the new method is not always the best, but in some
cases it is. It is also remarkable that empirically that the extended
Krylov-subspace approach requires so few iterations to \revised{build a} good
approximation of the solution \revised{set}; 
%\sout{manifold}
examination of the results in 
Appendix~\ref{sec:detailed numerical experiments} illustrates only 
a few instances where more than 20 iterations are required. We provide
analysis to explain why this should be so, although this may well be
tightened in future.

\revised{Our experiments suggest that, as anticipated, using an extended
Krylov subspace is most warranted when the required shift $\sigma$ lies
towards $\sigmamin$, while if $\sigma$ is large, a less-expensive, 
factorization-free standard Krylov space approach suffices. This then
opens the possibility of changing the frequency in which $A^kb$ and $A^{-\ell}$ 
terms are mixed.  Jagels and Reichel \cite{JageReic09:laa} have already
developed extended Krylov subspace methods that favour $A^kb$ terms, and
we are currently looking at more-general ideas along these lines
for our subproblems of interest.}

Software that implements the algorithms is now available 
\footnote{\revised{The {\tt TREK} package for the trust-region case,
and its counterpart {\tt NREK} for the the norm regularisation
case discussed in Section~\ref{regularization-section}, is
available with full documentation and test programs
see, \url{https://ralna.github.io/galahad_docs}.}}
as part of the {\sf GALAHAD} \cite{GoulOrbaToin03:toms} fortran library, 
and there are interfaces to other popular languages such as C, Python and Julia,
and tools such as Matlab. The codes are designed to solve not only
a single instance, but a sequence in which the regularization is tightened,
without recourse to expensive re-evaluations; the current subspace $V_k$
is extended if it is insufficient as the regularization changes.
The new packages will be rolled out as optional subproblem solvers
in other {\sf GALAHAD} packages in due course.

\section*{Acknowledgement}

\revised{The authors appreciate the helpful comments from two reviewers of 
this paper.}

\bibliographystyle{plain}
\bibliography{trek_paper}

\appendix
\section{Convergence analysis}
\label{sec:convergence analysis}
\renewcommand{\theequation}{A.\arabic{equation}}
\renewcommand{\thefigure}{A.\arabic{figure}}
The first $j$-th matrix products of the EKS algorithm 
(Algorithm~\ref{algorithm-eks}) generate an orthogonal basis for the
extended Krylov subspace $\calK_j(A,b)$. A natural question is then, how
closely do the elements of 
\revised{$\calS(A,b) = \mbox{span}\{\xs\}$} lie to \revised{$\calK_j(A,b)$}.
That is to say, for a given $\sigma$, is there an
entry in $\calK_j(A,b)$ close to $\xs$, and if so, how close?

Consider the symmetric positive-definite matrix $A$ with smallest and
largest eigenvalues $\lambda_1 = \lambda_1(A)$ and $\lambda_n = \lambda_n(A)$,
and consider the interval
$\calW = [\lambda_1,\lambda_n] \subset \mathbb{R_+^*}$.
In this section, we repeat for completeness, with minor variations, 
the theory developed by Knizhnerman and Simoncini \cite{KnizSimo11:NUMMAT} 
to prove that the convergence of 
the extended-Krylov-subspace method to build an approximate basis 
for the \revised{subspace} 
%\sout{solution manifold} 
$\calS(A,b)$, for all $\sigma \notin - \calW - ( -\delta_1,\delta_n) 
 = - (\lambda_1 - \delta_1,\lambda_n + \delta_n)$, 
and $\delta_1,\delta_n >0$, 
is linear at worst. The proof relies heavily
on 19th to mid 20th century complex analysis, and we refer the reader to 
\cite{KnizSimo11:NUMMAT} for the details; the technical buildup
may appear somewhat fearsome to the uninitiated.

Let $D$ denote the closed unit disk, let $\Psi :
\bar{\mathbb{C}}\setminus D \rightarrow \bar{\mathbb{C}}\setminus W$ be
the Riemann mapping for $\calW$ that preserves infinity and has positive
derivatives there and let $\Phi = \Psi^{-1}$ be its inverse. That is,
$\Psi = \Phi^{-1}$ and both $\Phi$ and $\Psi$ are holomorphic on their domains.  
For the specific set $\calW$, we have that
$\Phi(t)$ and $\Psi(t)$ are real \revised{for all} $t\notin \calW$,
and explicitly
$$
\Psi(t) = c +\frac{r}{2} \left(t+\frac{1}{t}\right) \text { and} 
$$ 
$$
\Phi(t) = 
\begin{cases}
    \bigfrac{t-c}{r} + \sqrt{\left(\bigfrac{t-c}{r}\right)^2-1} & 
 \text{ for } t \leq -\lambda_n\\
   \left(\bigfrac{t-c}{r} + \sqrt{\left(\bigfrac{t-c}{r}\right)^2-1}\right)^{-1} & 
 \text{ for } t \geq -\lambda_1
\end{cases}
$$
where 
$c = (\lambda_n+\lambda_1)/2$ and $r = (\lambda_n-\lambda_1)/2$.
We then define the function (a finite Blaschke product)
$$B(w) = w\frac{1-\overline{\Phi(0)}w}{\Phi(0)-w},$$
and note that this function satisfies $|B(w)| = 1$ for $|w|=1$ and 
\eqn{bw}{B(w)<1}
for $|w|<1$ (see, \cite{KnizSimo11:NUMMAT}).

Armed with these definitions, let $\phi_{2\ell}$ and $\phi_{2\ell + 1}$
be the (Takenaka–Malmquist) set of rational functions
%for the cyclic sequence of poles $0, \infty, 0, \infty, \ldots$.  
\eqn{tm}{\phi_{2\ell}(w) = B(w)^\ell \text{ and }
\phi_{2\ell +1}(w) = -\frac{\sqrt{1-|\Phi(0)|^{-2}}}{1-\Phi(0)^{-1}w} w B(w)^\ell.}
for $\ell \in\mathbb{N}$.
For these, we have
$$\frac{1}{2\pi}\int_{|w|=1}\overline{\phi_p(w)}\phi_m(w)|dw| = \delta_{m,p}$$
and
$$\max_{m\in \mathbb{N}}\max_{|w|=1} |\phi_m(w)| \leq c_{\phi}$$
for all $m,p\in\mathbb{N}$ and some finite constant $c_{\phi}$.

Now consider the Faber-Dzhrbashyan rational functions
$$M_m(z) = \frac{1}{2\pi i} \int_{\Gamma_{R}} \frac{\phi_m(\Phi(\xi))}{\xi-z}d\xi
 \text{ for } z \notin \overline{G_R} \text{ and } m\in\mathbb{N},$$
where
$$\Gamma_R = \{z \in \mathbb{C} \setminus W : |\Phi(z)| = R\}, \;\;
  G_R = W \cup \{z \in \mathbb{C}\setminus W: |\Phi(z)| < R\}$$
and $R>1$. In particular, $M_m$ is a Faber transformation of $\phi_m$,
and it then follows that
\eqn{mrat}{M_{2k}(z) = \frac{p_{2k}(z)}{z^{k}} \text{ and } 
 M_{2k+1}(z) = \frac{p_{2k+1}(z)}{z^{k+1}} }
where $p_m$ is a polynomial of degree $m$.
%$M_{2\ell}$ is a rational function of type $[2\ell/\ell]$ and $M_{2\ell+1}$ is
%a rational function of type $[2\ell + 1/\ell + 1]$ for $\ell \in \mathbb{N}$. 
%In addition, every finite pole of $M_{2\ell}$ and $M_{2\ell+1}$ is zero, and of
%multiplicity $\ell$ and $\ell + 1$, respectively. 
Furthermore, 
\eqn{Mm}{|M_m(z)|\leq c_M}
for some finite constant $c_M$  for all $m$ 
and all $z \in \calW$.

\begin{proposition}
Given $A$, $\Phi$, $\Psi$, $\phi_k$, $M_k$ and $\calW$ as above, we have that
\eqn{zmainv}{
(zI-A)^{-1} = \frac{1}{\Phi(z)\Psi^\prime(\Phi(z))}\sum_{k=0}^\infty 
 \phi_k(\Phi(z)^{-1})M_k(A) \text{ for every } z \notin \calW}
\end{proposition}

\noindent
\begin{proof}
%\bpr
In general, we have the identity \cite{Koch58}
\[\frac{\Psi^\prime (w)}{\Psi(w)-u} = \frac{1}{w} \sum_{k=0}^\infty 
 \overline{\phi_k\left(\frac{1}{\bar{w}}\right)}M_k(u).\]
The special (matrix-valued) case for which $u = A$ and 
\revised{$w = \Phi(z)I$, i.e.,} $\Psi(w) = z I$,
leads directly to \req{zmainv}.
%\epr
\end{proof}

%\noindent
%{\bf Just checking \ldots does $M_k(A)$ mean $U \diag(M_k(\lambda_i)) U^T$?} 

\begin{proposition}
\label{prop:convergence_rate}
Let $\sigma \in \mathbb{R} \setminus (-\lambda_n-\delta_n,-\lambda_1+\delta_1)$.
The EKS algorithm (Algorithm~\ref{algorithm-eks}) used to approximate 
$\xs = (A+\sigma I)^{-1}b$ converges at worst linearly with a rate 
\eqn{rate}{\mu = \sqrt{|E(\sigma)|} < 1,}
where $E: \mathbb{R}\setminus [-\lambda_n,-\lambda_1] \rightarrow \mathbb{R}$
is 
\[
    E(\sigma) \eqdef B\left(\Phi(-\sigma)^{-1}\right).
\]
Furthermore, if $\sigma \geq 0$, 
\eqn{mubnd}{\mu \leq \frac{\kappa^{1/4}-1}{\kappa^{1/4}+1},} 
where $\kappa = \lambda_n/\lambda_1$.
\end{proposition}

\begin{figure}
    \centering
    \includegraphics[width=0.49\linewidth]{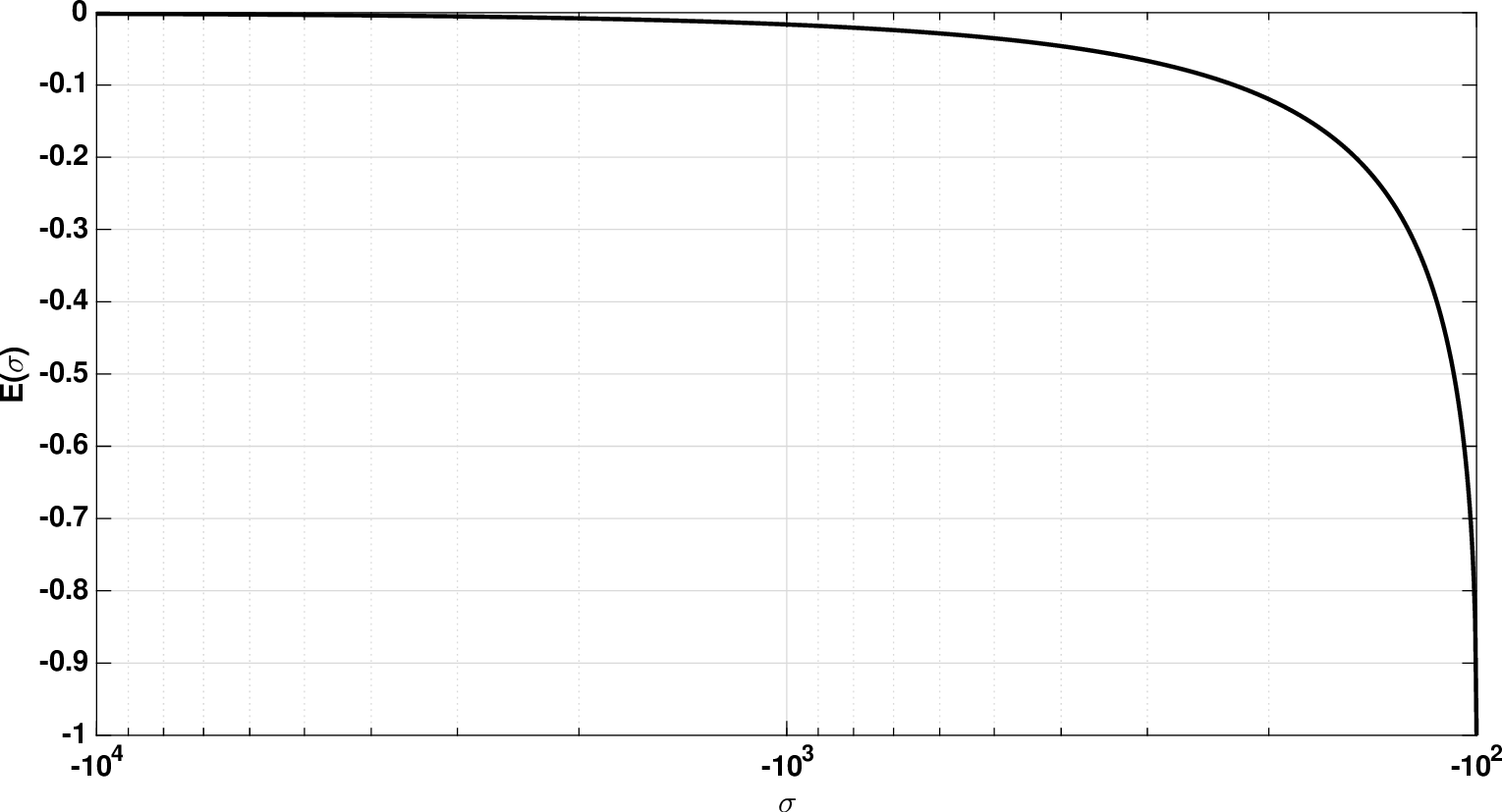}
    \includegraphics[width=0.49\linewidth]{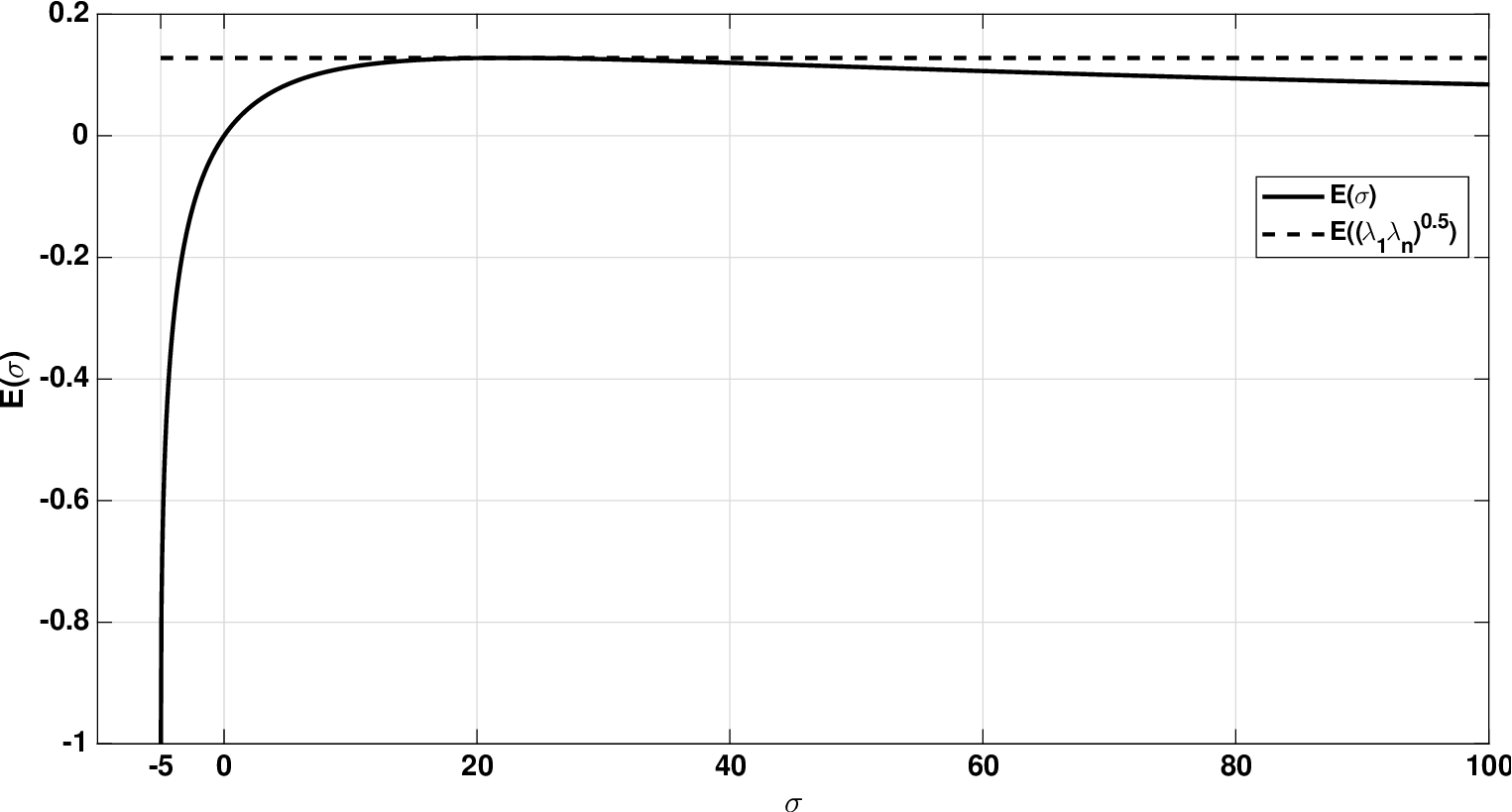}
    \caption{Upper bound on the convergence rate of extended Krylov method 
             to solve the shifted linear system with a shift $\sigma$. 
             The spectrum of the matrix is in $[5,100]$.}
    \label{fig:E}
\end{figure}

\noindent
\begin{proof}
%\bpr
Employing \req{zmainv} with $z = - \sigma$ immediately gives
\eqn{xs}{\begin{array}{rl}
x(\sigma) \ngap & = -((-\sigma) I - A)^{-1}b
  = \bigfrac{1}{\Phi(-\sigma)\Psi^\prime(\Phi(-\sigma))}
      \bigsum_{i=0}^\infty \phi_i(\Phi(-\sigma)^{-1})M_i(A)b \\
  & = \bigsum_{i=0}^\infty a_i M_i(A)b,
\end{array}}
where
\[a_i = \frac{\phi_i(\Phi(-\sigma)^{-1})}
 {\Phi(-\sigma)\Psi^\prime(\Phi(-\sigma))}.\]
Using \req{tm} with $w = -\sigma$ \revised{and studying the variation of $|a_i|$ as a function of $\sigma$}, we may bound
%Let us check how fast $|a_k|$ decays as $k$ grows. We have
\eqn{ak}{|a_i| \leq \frac{|\phi_i(\Phi(-\sigma)^{-1})|}{|\Phi(-\sigma)| 
 |\Psi^\prime(\Phi(-\sigma))|}
 \leq \frac{|B\left(\Phi(-\sigma)^{-1}\right)^{i/2}|}
 {|\Phi(-\sigma)| |\Psi^\prime(\Phi(-\sigma))|}
 \leq c_1 \mu^i,}
where
\[
c_1=\frac{1}{\Psi^\prime\left(\Phi(-\lambda_1-\delta_1)\right)}
\]
is a constant, and 
\eqn{mu}{\mu = \sqrt{|B\left(\Phi(-\sigma)^{-1}\right)|}
\equiv \sqrt{|E(\sigma)|} <1}
from \req{bw} as $|\Phi(-\sigma)^{-1}|<1$.

Since $M_i(A)b \in \calK_i(A,b)$ for $i=1,\ldots,2k$ 
by definition \req{eks} of $\calK_i(A,b)$ and from \req{mrat},
the estimate
\[x_j \eqdef \sum_{i=0}^{j} a_i M_i(A) b\]
is in $\calK_j(A,b)$ for $j = 1, \ldots 2k$.
But then we have
\[
\|x(\sigma)-x_j\| \leq \sum_{i=j+1}^\infty |a_i| \|M_i(A)\|\|b\| \leq c \mu^{j}.
\]
from \req{xs}, \req{Mm} and and \req{ak} respectively,
where $c = c_M c_1 \|b\| \mu / ( 1 - \mu)$ as $\mu < 1$.

A simple analysis, e.g., \cite[Prop.4.1]{KnizSimo11:NUMMAT}
shows that the function $E$ evolves as follows:
\vspace*{4mm}

\begin{tabular}{c|cccccccccccc}
   $\sigma$  & & $-\infty$ & &$-\lambda_n$ & & $-\lambda_1$ & & $0$ & &  
   $\sqrt{\lambda_1\lambda_n}$ & & $\infty$ \\
   \hline\\
     $E(\sigma)$ & & $0$ & $\searrow$ & 
    $-1$ & & $-1$ & $\nearrow$ & $0$ & $\nearrow$ & 
    $\left(\frac{\kappa^{1/4}-1}{\kappa^{1/4}+1}\right)^2$ & $\searrow$ & $0$ 
\end{tabular}

\vspace*{4mm}
\noindent
where $\kappa = \frac{\lambda_n}{\lambda_1}$; as an illustration, 
Figure~\ref{fig:E} 
depicts the function $E$ on $\mathbb{R}\setminus (-\lambda_n,-\lambda_1)$ 
with $\lambda_1 = 5$ and $\lambda_n = 100$.
Therefore, for symmetric positive-definite $A$ and a shift $\sigma\notin
[-\lambda_n,-\lambda_1]$, the convergence rate of EKS in approximating
the solution of the shifted linear system $(A+\sigma I)x=b$ is
$\sqrt{|E(\sigma)|}$, and if $\sigma\geq 0$, we have $E(\sigma) \leq
\left(\frac{\kappa^{1/4}-1}{\kappa^{1/4}+1}\right)^2$,
and thus \req{mubnd} follows from \req{mu}.
%\epr
\end{proof}

\begin{theorem}
The extended Krylov method to approximate a \revised{subspace containing}
%\sout{solution manifold}
\newline
$\{x(\sigma)\}_{\sigma\notin (-\lambda_n-\delta_n,-\lambda_1+\delta_1)}$ 
converges at worst linearly with a convergence rate 
$$\mu = 
 \max\left(\sqrt{|E(-\lambda_n-\delta_n)|},\sqrt{|E(-\lambda_1+\delta_1)|}, 
 \left(\frac{\kappa^{1/4}-1}{\kappa^{1/4}+1}\right)\right),$$
where $E$ is defined in Proposition~\ref{prop:convergence_rate}
and $\kappa = \lambda_n/\lambda_1$.
\end{theorem}

\noindent
\begin{proof}
%\bpr
The proof follows straightforwardly from 
Proposition~\ref{prop:convergence_rate}, and the discussion of the evolution
of $E(\sigma)$ given there.
%\epr
\end{proof}
%\clearpage

\section{Computing an orthogonal basis of $\calK_{\ell}(A,b)$}
\label{sec:orthogonal_basis}
\renewcommand{\theequation}{B.\arabic{equation}}
Given the $n$ by $n$  symmetric matrix $A$, the $n$-vector $b$,
we define a sequences of extended Krylov subspaces
\disp{\calK_{2k}(A,b) = \spanof\{A^{-k}b, \ldots, A^{k-1}b\}}
and
\disp{\calK_{2k+1}(A,b) = \spanof\{A^{-k}b, \ldots, A^{k}b\}
 =  \spanof\{\calK_{2k}(A,b), A^kb\}}
where $\calK_0(A,b) = \spanof\{b\}$, and thus
\disp{\calK_{2k}(A,b) = \spanof\{\calK_{2k-1}(A,b), A^{-k}b\}.}
We build a sequence of orthogonal basis for these so that
\disp{\calK_{2k}(A,b) = \spanof\{v_{-k}, \ldots, v_{k-1}\} \tim{and}
\calK_{2k+1}(A,b) = \spanof\{v_{-k}, \ldots, v_{k}\}}
with $\calK_0(A,b) = \{ b /\|b\|\}$.
Note that 
\disp{\arr{rl}{
\calK_{2k-1}(A,b) \ngap & = \spanof\{A^{-(k-1)}b, \ldots, A^{k-1}b\}
= \spanof\{v_{-(k-1)}, \ldots, v_{k-1}\}, \\
\calK_{2k-2}(A,b) \ngap & = \spanof\{A^{-(k-1)}b, \ldots, A^{k-2}b\}
= \spanof\{v_{-(k-1)}, \ldots, v_{k-2}\} \\
\tim{and} \calK_{2k-3}(A,b) \ngap & = \spanof\{A^{-(k-2)}b, \ldots, A^{k-2}b\}
= \spanof\{v_{-(k-2)}, \ldots, v_{k-2}\}}}
We find $v_{j}$ by orthogonalizing $A v_{-j}$ with respect to columns of $V_{2j}$,
in which case 
\eqn{Avmj}{A v_{-j} \in \calK_{2j+1},} 
and $v_{-(j+1)}$ by orthogonalizing $A^{-1}v_j$ with respect 
to the columns of $V_{2j+1}$, so that 
\eqn{Aivj}{A^{-1}v_j \in V_{2j+2}.}
By construction, it then follows that
\eqn{alpha}{\delta_k v_k = A v_{-k} - \sum_{j=-k}^{k-1} \alpha_j v_j
= A v_{-k} - \alpha_{-k} v_{-k} - \alpha_{k-1} v_{k-1} - \!\!\!
\sum_{j=-k+1}^{k-2} \alpha_j v_j}
for some normalizing constant $\delta_k$, which implies that
\disp{0 = \delta_k v_i^T v_{-k} = v_j^T A v_{-k} - \alpha_j,}
i.e.
\eqn{def-alpha}{\alpha_j = v_j^T A v_{-k} = v_{-k}^T A v_j}
for $j = -k,\ldots,k-1$. Consider first $j = -k+1,\ldots,-1$, and for these
\disp{\alpha_j \equiv \alpha_{-i} = v_{-i}^T A v_{-k} = v_{-k}^T A v_{-i}}
for $i = 1,\ldots,k-1$. But then, it follows from \req{Avmj} that
\disp{A v_{-i} \in \calK_{2i+1}(A,b) \subseteq \calK_{2k-1}(A,b),}
and thus $\alpha_{-i} = 0$ since $v_{-k}$ is orthogonal to 
$\calK_{2k-1}(A,b) =\spanof\{v_{-(k-1)}, \ldots, v_{k-1}\}$.
By contrast, consider $j = 0,\ldots,k-2$, and note that in this case
$v_j \in \calK_{2k-3}(A,b)$. But then
\disp{A v_j \in \spanof\{A^{-(k-3)}b, \ldots, A^{k-1}b\}
\subset \spanof\{A^{-(k-1)}b, \ldots, A^{k-1}b\} = \calK_{2k-1}(A,b)}
and thus as before $\alpha_{j} = 0$ since $v_{-k}$ is orthogonal to 
$\spanof\{v_{-(k-1)}, \ldots, v_{k-1}\}$. Thus \req{alpha} simplifies to
a short-term recurrence
\eqn{def-vk}{\delta_k v_k = A v_{-k} - \alpha_{-k} v_{-k} - \alpha_{k-1} v_{k-1}.}
Similarly, by construction, we have that
\eqn{beta}{\delta_{-(k+1)} v_{-(k+1)} = A^{-1} v_k - \sum_{j=-k}^k \beta_j v_j
= A^{-1} v_k - \beta_{-k} v_{-k} - \beta_k v_k -  \!\!\! 
\sum_{j=-k+1}^{k-1} \beta_j v_j}
for another normalizing constant $\delta_{-(k+1)}$, which implies that
\disp{0 = \delta_{-(k+1)} v_i^T v_{-(k+1)} = v_j^T A^{-1} v_{-k} - \beta_j,} 
i.e.
\eqn{def-beta}{\beta_j = v_j^T A^{-1} v_k = v_k^T A^{-1} v_j}
for $j = -k,\ldots,k$. Consider first $j = 0,\ldots,k-1$, and for these
\req{Aivj} implies that
\disp{A^{-1} v_j \in \calK_{2i+2}(A,b) \subseteq \calK_{2k}(A,b).}
But then $\beta_j = 0$ as $v_k$ does not lie in 
$\spanof\{v_{-k}, \ldots, v_{k-1}\}$. 
Whereas, if $j = -k+1,\ldots,-1$, 
\disp{\beta_j \equiv \beta_{-i} = v_{-i}^T A^{-1} v_k = v_k^T A^{-1} v_{-i}}
for $i = 1,\ldots,k-1$. Note that in this case $v_{-i} \in 
\calK_{2i}(A,b) \subseteq \calK_{2k-2}(A,b)$,
and then it follows that
\disp{A^{-1} v_{-i} \in \spanof\{A^{-k}b, \ldots, A^{k-3}b\}
\subset \spanof\{A^{-k}b, \ldots, A^{k-1}b\} = \calK_{2k}(A,b)}
Hence, once again, $\beta_j = 0$ as $v_k$ is orthogonal to 
$\spanof\{v_{-k}, \ldots, v_{k-1}\}$. Thus, as before, \req{beta} reduces
to a short-term recurrence
\eqn{def-vmk}{\delta_{-(k+1)} v_{-(k+1)} = 
 A^{-1} v_k - \beta_{-k} v_{-k} - \beta_k v_k.}
The identities \req{def-alpha}, \req{def-vk}, \req{def-beta}
and \req{def-vmk} then lead to Algorithm~\ref{algorithm-eks}.

\section{Modified norm trust region subproblem solver}
\label{sec:modified norm subproblem solver}
\renewcommand{\theequation}{C.\arabic{equation}}
\renewcommand{\thealgocf}{C.\arabic{algocf}}
\renewcommand{\thefigure}{C.\arabic{figure}}
\renewcommand{\thetable}{C.\arabic{table}}
\setcounter{algocf}{0}
Here we provide a detailed description of the modifications to the {\tt TREK}
algorithm (Algorithm~\ref{algorithm-trek}) that are necessary when using
the trust region norm for which $\|x\|_S^2 = x^T S x$, with $S = \PP^T \PP$,
as discussed in Section~\ref{other-elliptical-norms-section}.
As we stated there, it suffices to apply Algorithm~\ref{algorithm-trek}
with the matrix $\Ap = \PP^{-T} A \PP^{-1}$ and vector $\bp = \PP^{-T} b$, 
and to recover the solution via $x_* = \PP^{-1} \xpstar$. Since 
$S^{-1} = \PP^{-1} \PP^{-T}$, this is formal\revised{ized} as Algorithm~\ref{algorithm-trekp}
\vpageref{algorithm-trekp}.

\begin{algorithm2e}[!ht]
\caption{The {\tt TREK} algorithm to solve the trust-region subproblem
\req{trp} \label{algorithm-trekp}}
\KwInput{symmetric $A \in \Re^{n\times n}$, $b \in \Re^n$,
 symmetric, positive-definite $S \in \Re^{n\times n}$, $\Delta > 0$,
 a stopping tolerance $\epsilon > 0$  and an iteration bound $m\geq 1$.}
\KwOutput{$x_* = \approxargmininx{} \half x^T A x - b^T x 
 \tim{subject to} \|\PP x\| \leq \Delta.$}
  \setstretch{0.8}
  \vspace{2mm}
  {\bf Algorithm:}\\
  $x \define A^{-1} b$, $w \define \PP x$\;
  \If{$\|w\| \leq \Delta$}{
  {\bf exit} with the interior solution $x_* = x$ and 
    shift $\sigma_* = 0$\;
  }
  $\delta_0 \define \|\PP^{-T} b\|$, 
  $v_{0} \define \PP^{-T} b/\delta_0$, 
  $u \define w / \delta_0$\;
  $\beta_{0} \define u^T v_0$\;
  $u \redef u - \beta_{0} v_{0}$\;
  $\delta_{-1} \define \|u\|$\;
  \For{$k = 1, 2, \ldots , m$}{
    \uIf{$\delta_{-k} > 0$}{
     $v_{-k} \define u/\delta_{-k}$\;
     $u \define \PP^{-T} A \PP^{-1} v_{-k}$\;
     $\alpha_{k-1} \define u^T v_{k-1}$,
     $u \redef u - \alpha_{k-1} v_{k-1}$\;
     $\alpha_{-k} \define u^T v_{-k}$,
     $u \redef u - \alpha_{-k} v_{-k}$\;
     $\delta_{k} \define \|u\|$\;
    }
    \uIf{$k>1$}{update $P_{2k-1}$ from $P_{2k-2}$ using \req{podd} and form
     $b_{2k-1} = (b^T_{2k-2},0)^T$\;}
    \Else{initialize $P_1$ from \req{pstart} and set $b_1 = \delta_0$\;}
    $y_{2k-1} \define \argmin \half y^TP_{2k-1} y - b_{2k-1}^Ty \; \mbox{s.t.} \;
      \|y\| \leq \Delta$ and its optimal shift $\sigma_{2k-1}$\;
    compute $\|r_{2k-1}\|$ from \req{rkp}\;
    \If{$\|r_{2k-1}\| \leq \epsilon$}{
    {\bf exit} with the interior solution $x_* = V_{2k-1} y_{2k-1}$, shift
      $\sigma_* = \sigma_{2k-1}$ and $k_{\mbox{\scriptsize stop}} = k$\;}
    \uIf{$\delta_{k} > 0$}{
     $v_{k} \define u/\delta_{k}$, 
     $u \define \PP A^{-1} \PP^T v_{k}$\;
     $\beta_{-k} \define u^T v_{-k}$,
     $u \redef u - \beta_{-k} v_{-k}$\;
     $\beta_{k} \define u^T v_{k}$,
     $u \redef u - \beta_{k} v_{k}$\;
    } 
    \If{$k<m$}{\uIf{$\delta_{k} > 0$}{$\delta_{-k-1} \define \|u\|$\;}
     update $P_{2k}$ from $P_{2k-1}$ using \req{peven} and form
      $b_{2k} = (b^T_{2k-1},0)^T$\;
     $y_{2k} \define \argmin \half y^TP_{2k} y - b_{2k}^Ty \; \mbox{s.t.} \;
       \|y\| \leq \Delta$ together with its optimal shift $\sigma_{2k}$\;
     compute $\|r_{2k}\|$ from \req{rk}\;
     \If{$\|r_{2k}\| \leq \epsilon$}{
     {\bf exit} with the interior solution $x_* = V_{2k} y_{2k}$, shift
       $\sigma_* = \sigma_{2k}$ and $k_{\mbox{\scriptsize stop}} = k$\;
     }
  }
 }
\end{algorithm2e}

Under the transformation  $w_k = \PP^{-1} v_k$ i.e., $v_k = \PP w_k$, and
$u = \PP v$, i.e., $v = \PP^{-1}u$ (and then replacing $v_k \rightarrow w_k$ 
and $u \rightarrow v$ so that the resulting changes may be easily 
compared with Algorithm~\ref{algorithm-trek}),
and introducing temporary workspace vectors
$z$ and $s$, we may rewrite Algorithm~\ref{algorithm-trekp} as
Algorithm~\ref{algorithm-treks} \vpageref{algorithm-treks}.

\begin{algorithm2e}[!ht]
\caption{The {\tt TREK} algorithm to solve the trust-region subproblem
\req{trss} \label{algorithm-treks}}
\KwInput{symmetric $A \in \Re^{n\times n}$, $b \in \Re^n$, 
 symmetric, positive-definite $S \in \Re^{n\times n}$, $\Delta > 0$,
 a stopping tolerance $\epsilon > 0$  and an iteration bound $m\geq 1$.}
\KwOutput{$x_* = \approxargmininx{} \half x^T A x - b^T x 
 \tim{subject to} \|x\|_S \leq \Delta.$}
  \setstretch{0.8}
  \vspace{2mm}
  {\bf Algorithm:}\\
  $x \define A^{-1} b$, 
  $u \define S x$\;
  \If{$\sqrt{u^T x} \leq \Delta$}{
  {\bf exit} with the interior solution $x_* = x$ and 
    shift $\sigma_* = 0$\;
  }
  $w \define S^{-1} b$, 
  $\delta_0 \define \sqrt{b^T w}$\;
  $v_0 \define w / \delta_0$, 
  $w \define u / \delta_0$,
  $u \define x / \delta_0$\;
  $\beta_0 \define v_{0}^T S u$, 
  $u \redef u - \beta_0 v_0$\;
  $w = S u$, 
  $\delta_{-1} \define \sqrt{u^T w}$\;
  \For{$k = 1, 2, \ldots , m$}{
    \uIf{$\delta_{-k} > 0$}{
     $v_{-k} \define u/\delta_{-k}$,
     $w \define A v_{-k}$, 
     $u \define S^{-1} w$\;
     $\alpha_{-k} \define v_{-k}^T S u$,
     $u \redef u - \alpha_{-k} v_{-k}$\;
     $\alpha_{k-1} \define v_{k-1}^T S u$,
     $u \redef u - \alpha_{k-1} v_{k-1}$\;
     $w = S u$, 
     $\delta_{k} \define \sqrt{u^T w}$\;
    }
    \uIf{$k>1$}{update $P_{2k-1}$ from $P_{2k-2}$ using \req{podd} and form
     $b_{2k-1} = (b^T_{2k-2},0)^T$\;}
    \Else{initialize $P_1$ from \req{pstart} and set $b_1 = \delta_0$\;}
    $y_{2k-1} \define \argmin \half y^TP_{2k-1} y - b_{2k-1}^Ty \; \mbox{s.t.} \;
      \|y\| \leq \Delta$ and its optimal shift $\sigma_{2k-1}$\;
    compute $\|r_{2k-1}\|$ from \req{rkp}\;
    \If{$\|r_{2k-1}\| \leq \epsilon$}{
    {\bf exit} with the interior solution $x_* = V_{2k-1} y_{2k-1}$, shift
      $\sigma_* = \sigma_{2k-1}$ and $k_{\mbox{\scriptsize stop}} = k$\;}
    \uIf{$\delta_{k} > 0$}{
     $v_{k} \define u/\delta_{k}$,
     $w \define  S v_{k}$,
     $u \define  A^{-1} w$\;
     $\beta_{-k} \define v_{-k}^T S u $,
     $u \redef u - \beta_{-k} v_{-k}$\;
     $\beta_{k} \define v_{k}^T S u$,
     $u \redef u - \beta_{k} v_{k}$\;
    }
    \If{$k<m$}{\uIf{$\delta_{k} > 0$}{$w = S u$, 
     $\delta_{-k-1} \define \sqrt{u^T w}$\;
     }
     update $P_{2k}$ from $P_{2k-1}$ from \req{P} and form
      $b_{2k} = (b^T_{2k-1},0)^T$\;
     $y_{2k} \define \argmin \half y^TP_{2k} y - b_{2k}^Ty \; \mbox{s.t.} \;
       \|y\| \leq \Delta$ together with its optimal shift $\sigma_{2k}$\;
     compute $\|r_{2k}\|$ from \req{rk}\;
     \If{$\|r_{2k}\| \leq \epsilon$}{
     {\bf exit} with the interior solution $x_* = V_{2k} y_{2k}$, shift
       $\sigma_* = \sigma_{2k}$ and $k_{\mbox{\scriptsize stop}} = k$\;
     }
  }
 }
\end{algorithm2e}

%\noindent
While it might appear that Algorithm~\ref{algorithm-treks} requires multiple
products of $S$ with selected vectors, as needed to calculate the 
terms $v_i^T S u$ that define $\alpha_i$ and $\beta_i$, 
fortunately this is not the case.
There are two ways to avoid this. In the first, we maintain appropriate
vectors $q_i \eqdef S v_i$. We \revised{formalize} this as Algorithm~\ref{algorithm-treks1}
\vpageref{algorithm-treks1}.

\begin{algorithm2e}[!ht]
\caption{The {\tt TREK} algorithm to solve the trust-region subproblem
\req{trss} \label{algorithm-treks1}}
\KwInput{symmetric $A \in \Re^{n\times n}$, $b \in \Re^n$, 
 symmetric, positive-definite $S \in \Re^{n\times n}$, $\Delta > 0$,
 a stopping tolerance $\epsilon > 0$  and an iteration bound $m\geq 1$.}
\KwOutput{$x_* = \approxargmininx{} \half x^T A x - b^T x 
 \tim{subject to} \|x\|_S \leq \Delta.$}
  \setstretch{0.8}
  \vspace{2mm}
%  \small
  {\bf Algorithm:}\\
  $x \define A^{-1} b$, 
  $u \define S x$\;
  \If{$\sqrt{u^T x} \leq \Delta$}{
  {\bf exit} with the interior solution $x_* = x$ and 
    shift $\sigma_* = 0$\;
  }
  $w \define S^{-1} b$, 
  $\delta_0 \define \sqrt{b^T w}$\;
  $v_0 \define w / \delta_0$, 
  $q_0 \define b/\delta_0$, 
  $w \define u / \delta_0$,
  $u \define x / \delta_0$\;
  $\beta_0 \define u^T q_0$, 
  $u \redef u - \beta_0 v_0$\;
  $w = S u$, 
  $\delta_{-1} \define \sqrt{u^T w}$\;
  \For{$k = 1, 2, \ldots , m$}{
    \uIf{$\delta_{-k} > 0$}{
     $v_{-k} \define u/\delta_{-k}$,
     $q_{-k} \define w/\delta_{-k}$\;
     $w \define A v_{-k}$, $u \define S^{-1} w$\;
     $\alpha_{-k} \define u^T q_{-k}$,
     $u \redef u - \alpha_{-k} v_{-k}$\;
     $\alpha_{k-1} \define u^T q_{k-1}$,
     $u \redef u - \alpha_{k-1} v_{k-1}$\;
     $w = S u$, 
     $\delta_{k} \define \sqrt{u^T w}$\;
    }
    \uIf{$k>1$}{update $P_{2k-1}$ from $P_{2k-2}$ using \req{podd} and form
     $b_{2k-1} = (b^T_{2k-2},0)^T$\;}
    \Else{initialize $P_1$ from \req{pstart} and set $b_1 = \delta_0$\;}
    $y_{2k-1} \define \argmin \half y^TP_{2k-1} y - b_{2k-1}^Ty \; \mbox{s.t.} \;
      \|y\| \leq \Delta$ and its optimal shift $\sigma_{2k-1}$\;
    compute $\|r_{2k-1}\|$ from \req{rkp}\;
    \If{$\|r_{2k-1}\| \leq \epsilon$}{
    {\bf exit} with the interior solution $x_* = V_{2k-1} y_{2k-1}$, shift
      $\sigma_* = \sigma_{2k-1}$ and $k_{\mbox{\scriptsize stop}} = k$\;}
    \uIf{$\delta_{k} > 0$}{
     $v_{k} \define u/\delta_{k}$,
     $q_{k} \define w/\delta_{k}$\;
     $u \define  A^{-1} q_k$\;
     $\beta_{-k} \define u^T q_{-k}$,
     $u \redef u - \beta_{-k} v_{-k}$\;
     $\beta_{k} \define u^T q_{k}$,
     $u \redef u - \beta_{k} v_{k}$\;
    }
    \If{$k<m$}{\uIf{$\delta_{k} > 0$}{$w = S u$, 
     $\delta_{-k-1} \define \sqrt{u^T w}$\;
    }
     update $P_{2k}$ from $P_{2k-1}$ using \req{peven} and form
      $b_{2k} = (b^T_{2k-1},0)^T$\;
     $y_{2k} \define \argmin \half y^TP_{2k} y - b_{2k}^Ty \; \mbox{s.t.} \;
       \|y\| \leq \Delta$ together with its optimal shift $\sigma_{2k}$\;
     compute $\|r_{2k}\|$ from \req{rk}\;
     \If{$\|r_{2k}\| \leq \epsilon$}{
     {\bf exit} with the interior solution $x_* = V_{2k} y_{2k}$, shift
       $\sigma_* = \sigma_{2k}$ and $k_{\mbox{\scriptsize stop}} = k$\;
     }
  }
 }
\end{algorithm2e}
%\noindent
In practice, use two arrays $q_+$ and $q_-$, and store and successively 
overwrite $q_+ = q_k$ for $k \geq 0$ and $q_- = q_k$ for $k < 0$ with
$k = 1, 2, \ldots , m$.

The second way to avoid multiple products with $S$ is to use $q_i$ as before, 
but now to update $w \eqdef S u$ as it appears. We record this variant as 
Algorithm~\ref{algorithm-treks2} \vpageref{algorithm-treks2}.
Once again, in practice, use two arrays $q_+$ and 
$q_-$, and store and successively overwrite $q_+ = q_k$ for $k \geq 0$ and 
$q_- = q_k$ for $k < 0$.

Computationally, the main difference is that Algorithm~\ref{algorithm-treks1} 
requires an additional $S$ product\revised{s} per cycle, while 
Algorithm~\ref{algorithm-treks2} has an additional four vector subtractions.

\clearpage

%\section*{Appendix 2}
\section{Detailed numerical experiments}
\label{sec:detailed numerical experiments}
\renewcommand{\thetable}{D.\arabic{table}}
Here we summarize the results from our numerical experiments that lead to
Figure~\ref{figure-times}. As we said, we consider the set of \revised{93} 
unconstrained 
test problems in the current release (2025-09-09) of the \cutest\ optimization 
test examples \cite{GoulOrbaToin15:coap} that have 1000 or more variables, 
evaluate the gradient and  Hessian to provide $-b$ and $A$, 
\revised{and for 
each two sets of radii (a) $\Delta_0 = 10$ and then $\Delta = 1$ and $0.1$,
and (b) $\Delta_0 = 1$ and then $0.1$ are used.\footnote{In two exceptional 
cases, for problems {\tt ARWHEAD} and {\tt EG2}, both $\Delta = 10$ and $1$ 
give interior solutions, and for these we pick $\Delta_0 = 10$ and then 
$\Delta = 0.1$ and $0.01$ in (a).}
}
We apply the \galahad\ solvers 
{\tt TREK} (Algorithm~\ref{algorithm-trek}),
{\tt TRS} \cite{GoulRobiThor10:mpc}
and {\tt GLTR} \cite{GoulLuciRomaToin99:siopt} 
in each case \revised{for both (a) and (b) regimes; when $\Delta$ is smaller 
than its initial value, we use data accumulated from the solve with the 
initial radius case to ``warm start'' the second (and third) solves 
as each solver provides for this.}

In Table~\ref{table-a1} we report the name of each test problem, the size of
the instance chosen, $n$ (many \cutest\ examples have multiple dimensions),
the radius value, $\Delta$ and the optimal objective function found.
For {\tt TREK}, we give the terminal iteration count, $k_{\mbox{\scriptsize stop}}$,
which is the number of solves with $A$ required, and the clock time required.
For {\tt TRS}, we report the number of factorizations and clock time required,
while for {\tt GLTR}, we state the number of vector products involving $A$
and again the time. Both of the $\Delta$ regimes (a) and (b) are included
for each solver. The symbol $\dag$ indicates an instance 
that the iteration bound \revised{$m = 300$} was reached, while
and $\ddag$ indicates that a factorization failed through singularity.

{\tiny
%{\fontsize{1}{1}\selectfont
\setlength{\tabcolsep}{3pt}
\begin{longtable}{|lrrr|rr|rr|rr|rr|rr|rr|}
\caption{\label{table-a1}{Subproblem solves for CUTEst examples: Complete results}}
 \\
\hline
           &        &     &                 & \multicolumn{4}{c|}{TREK} & \multicolumn{4}{c|}{TRS} & \multicolumn{4}{c|}{GLTR} \\ \cline{5-16}
%           &        &     &                 & \multicolumn{2}{c|}{$\Delta_0=10$} & \multicolumn{2}{c|}{$\Delta_0=1$} & \multicolumn{2}{c|}{$\Delta_0=10$} & \multicolumn{2}{c|}{$\Delta_0=1$} & \multicolumn{2}{c|}{$\Delta_0=10$} & \multicolumn{2}{c|}{$\Delta_0=1$} \\
           &        &     &                 & \multicolumn{2}{c|}{(a)} & \multicolumn{2}{c|}{(b)} & \multicolumn{2}{c|}{(a)} & \multicolumn{2}{c|}{(b)} & \multicolumn{2}{c|}{(a)} & \multicolumn{2}{c|}{(b)} \\
name       &    $n$   & $\Delta$ & $f(x_*)$ &   $k_{\mbox{\tiny s}}$ &   time &  $k_{\mbox{\tiny s}}$ &   time  & fac &   time & fac &   time &  pr &   time &  pr &   time \\
\hline
\endfirsthead
\caption{{Subproblem solves for CUTEst examples: Complete results (continued)}} \\
\hline
           &        &     &                 & \multicolumn{4}{c|}{TREK} & \multicolumn{4}{c|}{TRS} & \multicolumn{4}{c|}{GLTR} \\ \cline{5-16}
           &        &     &                 & \multicolumn{2}{c|}{(a)} & \multicolumn{2}{c|}{(b)} & \multicolumn{2}{c|}{(a)} & \multicolumn{2}{c|}{(b)} & \multicolumn{2}{c|}{(a)} & \multicolumn{2}{c|}{(b)} \\
name       &    $n$   & $\Delta$ & $f(x_*)$ &   $k_{\mbox{\tiny s}}$ &   time &  $k_{\mbox{\tiny s}}$ &   time  & fac &   time & fac &   time &  pr &   time &  pr &   time \\
\hline
\endhead
\hline
\endfoot
\hline
\endlastfoot
ARWHEAD    &   5000 & 10  & -9.99800000E+03 &    0       &     0.01 &                             2        &     0.00 &             &          &      2 &     0.00 &        &          \\
ARWHEAD    &   5000 & .1  & -3.59936000E+03 &    1       &     0.00 &     0       &     0.01 &    4        &     0.00 &           2 &     0.00 &      1 &     0.00 &      2 &     0.00 \\
ARWHEAD    &   5000 & .01 & -3.95930600E+02 &    1       &     0.00 &     1       &     0.00 &    6        &     0.00 &           4 &     0.00 &      1 &     0.00 &      1 &     0.00 \\
BA-L16LS   &  66462 & 10  & -3.46265942E+11 &   12       &     0.91 &             &          &    9        &     1.28 &             &          &     75 &     2.11 &        &          \\
BA-L16LS   &  66462 & 1   & -4.00647969E+09 &   12       &     0.00 &    10       &     0.82 &   13        &     0.79 &           6 &     1.06 &     37 &     1.08 &     69 &     1.96 \\
BA-L16LS   &  66462 & .1  & -2.71798703E+08 &   12       &     0.00 &    10       &     0.00 &   19        &     0.98 &          12 &     1.02 &     37 &     1.09 &     34 &     1.01 \\
BA-L21LS   &  34134 & 10  & -1.09802044E+09 &  284       &    11.60 &             &          &   15        &     0.70 &             &          &   1999 &    22.79 &        &          \\
BA-L21LS   &  34134 & 1   & -1.43350953E+08 &  284       &     0.15 &   196       &     3.34 &   21        &     0.39 &           8 &     0.43 &    999 &    11.38 &   1129 &    13.46 \\
\hline
\end{longtable}
}

%\clearpage
%\appendix 

\end{document}